\documentclass[11pt,reqno,oneside]{amsart}
 \usepackage{verbatim,amsmath,amssymb,cite,epsfig,xspace}

\numberwithin{equation}{section}

\theoremstyle{plain}
\newtheorem{theorem}{Theorem}[section]
\newtheorem{lemma}{Lemma}[section]
\newtheorem{corollary}{Corollary}[section]
\theoremstyle{definition}

\theoremstyle{remark}
\newtheorem{remark}{Remark}[section]

\newcommand{\Real}{\mathbb R}

\newcommand{\pb}{\mathbf \Psi}
\newcommand{\Pb}{\mathbf \Phi}
\newcommand{\ab}{\mathbf a_0}
\newcommand{\az}{\mathbf a}
\newcommand{\fb}{\mathbf f}
\newcommand{\gb}{\mathbf g}
\newcommand{\hb}{\mathbf h}

\newcommand{\rb}{\mathbf r}
\newcommand{\Rb}{\mathbf R}
\newcommand{\sbf}{\mathbf s}
\newcommand{\yb}{\mathbf y}
\newcommand{\zb}{\mathbf z}

\newcommand{\E}{{\mathcal E}}
\newcommand{\Po}{{\mathcal P}}
\newcommand{\Su}{{\mathcal S}}

\def\maxnorm#1{||#1||_\infty}

\newcommand{\half}{\frac{1}{2}}

\begin{document}

\title[Analysis of a  Force-Based Quasicontinuum Approximation]
{Analysis of a Force-Based Quasicontinuum Approximation}
\author{Matthew Dobson}
\author{Mitchell Luskin}

\address{Matthew Dobson\\
School of Mathematics \\
University of Minnesota \\
206 Church Street SE \\
Minneapolis, MN 55455 \\
U.S.A.}
\email{dobson@math.umn.edu}

\address{Mitchell Luskin \\
School of Mathematics \\
University of Minnesota \\
206 Church Street SE \\
Minneapolis, MN 55455 \\
U.S.A.}
\email{luskin@umn.edu}

\thanks{This work was supported in part by
 DMS-0304326
  and by the Minnesota Supercomputer Institute.
  This work is also based on
work supported by the Department of Energy under Award Number
DE-FG02-05ER25706.
}



\date{\today}


\begin{abstract}
We analyze a force-based quasicontinuum approximation to a
one-dimensional system of atoms that interact by a classical
atomistic potential. This force-based quasicontinuum approximation
can be derived as the modification of an energy-based
quasicontinuum approximation by the addition of nonconservative
forces to correct nonphysical ``ghost'' forces that occur
in the
atomistic to continuum interface.  The algorithmic simplicity
and improved accuracy of the
force-based quasicontinuum approximation has made it
popular for large-scale quasicontinuum computations.

We prove that the force-based quasicontinuum equations have
a unique solution when the magnitude of the external forces satisfy
explicit bounds.   For Lennard-Jones next-nearest-neighbor
interactions, we show that unique solutions exist
for external forces that extend the system nearly to its tensile limit.

We give an analysis of the convergence of the ghost force iteration method
to solve the equilibrium equations for the force-based quasicontinuum approximation.
We show that the ghost force iteration is a contraction and give an analysis for its
convergence rate.
\end{abstract}

\maketitle
{
\thispagestyle{empty}

\section{Introduction}

The local lattice structure for minimum energy configurations
of atomistic systems subject to external
forces is usually slowly varying except near defects such as
dislocations~\cite{tadmor_miller_qc_overview}.
Quasicontinuum methods efficiently approximate
these multiscale features by maintaining
atomistic degrees of freedom near defects and coarse-graining
the atomistic degrees of freedom in regions where the local lattice
structure is nearly uniform through the introduction of
representative atoms~
\cite{tadmor_miller_qc_overview,tadmor_qc_first,knaportiz}.
The efficiency of quasicontinuum methods has
allowed the simulation of more complex
problems than can be computed using a completely atomistic
model~\cite{shenoy_gf}.

Many quasicontinuum methods have been
proposed~\cite{tadmor_qc_first, knaportiz,ezhang06,
tadmor_miller_qc_overview, rodney_gf, miller_indent, ortnersuli, jacobsen04},
and each version gives a
different quasicontinuum approximation of the atomistic system.
A force-based quasicontinuum approximation has been proposed that
modifies an energy-based quasicontinuum approximation
by the addition of nonconservative forces to correct nonphysical
``ghost'' forces that occur in the atomistic to continuum interface
~\cite{tadmor_miller_qc_overview,rodney_gf,shenoy_gf,knaportiz}.
The force-based quasicontinuum approximation has been popular
for large-scale computations
since the improved accuracy is obtained with no additional
computational work simply by computing the
force on each representative atom with either
an atomistic algorithm or with
a continuum finite element algorithm
~\cite{tadmor_miller_qc_overview,rodney_gf,shenoy_gf,knaportiz,PrudhommeBaumanOden:2005}.

Adaptive mesh and error control have been successfully used
with the force-based quasicontinuum approximation
to efficiently choose representative atoms
~\cite{shenoy_gf,knaportiz,PrudhommeBaumanOden:2005,OdenPrudhommeRomkesBauman:2005}.
The number of representative atoms surrounding defects that need to be
modeled atomistically (the core of the defect) can be determined by the 
error tolerance, and the mesh in the continuum region surrounding the core 
can be coarsened beyond the atomistic-continuum interface wherever the 
deformation gradient varies slowly.  For simplicity,
reported implementations have not coarsened within the
cut-off radius of atomistic representative atoms, but it is possible
to coarsen immediately beyond the atomistic-continuum interface by
interpolating between continuum representative atoms.

In Section~\ref{model},
we give a derivation following~\cite{tadmor_miller_qc_overview,rodney_gf,shenoy_gf,knaportiz}
of several quasicontinuum
approximations leading to the derivation of
the force-based quasicontinuum approximation.
In Section~\ref{integ}, we reformulate the equilibrium equations
as a balance of forces conjugate to the distances between representative
atoms, rather than as a balance of forces conjugate to the
positions of the representative atoms.  Our derivation
and reformulation gives the
mathematical structure that is used in
our analysis.

In Section~\ref{exist}, we prove that the force-based
quasicontinuum equations have a unique solution under suitable
restrictions on the loads.  In the case of Lennard-Jones
next-nearest-neighbor interactions, we determine bounds
for the magnitude of the loads for which
unique solutions exist and find that the allowable loads
extend quite close to the tensile limit.

In Section~\ref{iterza}, we give an analysis of the convergence
of the ghost force iteration method that has been most commonly used
to solve the equilibrium equations for the
force-based quasicontinuum approximation
~\cite{tadmor_miller_qc_overview,rodney_gf,shenoy_gf,knaportiz}.
We prove that the ghost force iteration is a contraction and give
a bound for its
convergence rate.  We show that our bound for the convergence rate
gives a high convergence rate
when applied to the Lennard-Jones model subject to moderate
external forces.

Mathematical analyses of energy-based versions of the
quasicontinuum approximation that do not include ghost force
corrections have been given in
~\cite{pinglin03,pinglin05,legollqc05,ortnersuli,ortnersuli2,minge04,minge05},
and a simplified version of our analysis can be used to prove
the existence of solutions to these energy-based quasicontinuum approximations.
 We show, though, that the ghost forces are
nonconservative forces, so they cannot be derived from an energy.
Thus, the force-based quasicontinuum approximation cannot be
completely analyzed by energy methods.

We refer to \cite{lions} for a review of current progress on
the mathematical analysis of atomistic to continuum models for
solids and to \cite{tadmor_miller_qc_overview} for an introduction
and overview of the quasicontinuum approximation.

\section{Quasicontinuum Approximations}
\label{model}

In this section, we describe a sequence of
one-dimensional coarse-grained approximations
of a chain of atoms with nearest-neighbor and
next-nearest-neighbor interactions
given by a classical two-body potential, $\phi(r).$  We assume that
the atomistic potential
$\phi(r)$ is defined for all $r>0.$

We begin with the atomistic model, which has degrees of freedom for
all atomic positions
and computes a total internal energy directly from pairwise interactions.  From there,
we examine the constrained atomistic approximation and
the local quasicontinuum approximation which both
decrease the degrees of freedom by interpolating atomic positions between
representative atoms.
We then introduce an energy-based and a force-based quasicontinuum
approximation which span atomistic and continuum scales by combining
atomistic regions where the atoms directly interact according to the
atomistic model and continuum regions where the atoms interact according
to the local quasicontinuum approximation.  We observe that the
energy-based quasicontinuum approximation gives nonphysical
ghost forces near the atomistic to continuum interface that are
corrected by the force-based quasicontinuum approximation.

\subsection{The Atomistic Model}
\begin{figure}[hb]
\includegraphics{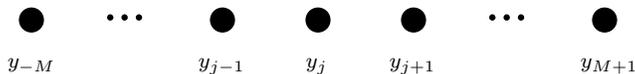}
\caption{Atomistic chain with atoms labeled by their position.}
\end{figure}
We denote the positions of the atoms by
$y_{i}$ for $i=-M,\dots,M+1,$ where
$
y_{i}<y_{i+1}.
$
The total energy for the atomistic system with nearest-neighbor and
next-nearest-neighbor interactions given by the classical two-body
potential $\phi(r)$ is
\begin{equation*}
\E^a(\yb) = \sum_{i=-M}^{M} \left[
\phi(y_{i+1} - y_{i})+\phi(y_{i+2}-y_{i})
\right],
\end{equation*}
where $\yb=(y_{-M},\dots,y_{M+1})\in\Real^{2M+2}$
and where the boundary terms
$\phi(y_{i} - y_{j})$ above and in the following
should be understood to be zero
for $i\notin \{-M,\dots,M+1\}$ or $j\notin \{-M,\dots,M+1\}.$
We can also express the total energy in terms of energies associated with
each atom as
\begin{equation}\label{atom2}
\E^a(\yb) = \sum_{i=-M}^{M+1} \E_i^a(\yb)
\end{equation}
with
\begin{equation} \label{atompart}
\E_i^a(\yb) =
\half\Big[\phi(y_{i+1}-y_{i})
+\phi(y_{i +2}-y_{i})
+ \phi(y_{i}-y_{i-1})
+ \phi(y_{i}-y_{i-2})\Big].
\end{equation}
We then have that the force on the atom at position $y_{i}$
is given by
\begin{equation} \label{atomforce}
\begin{split}
F_{i}^a(\yb) &= -\frac{\partial}{\partial y_{i}}
\Big[\phi(y_{i +1}-y_{i })
+ \phi(y_{i +2}-y_{i})
+ \phi(y_{i }-y_{i-1})
+ \phi(y_{i }-y_{i -2})\Big] \\
&= \left[\eta(r_{i}) + \eta(r_{i}+r_{i+1})\right]
 - \left[\eta(r_{i-1}) + \eta(r_{i-1}+r_{i-2})\right],
\end{split}
\end{equation}
where $\eta(r) = \phi'(r)$ and $r_i=y_{i+1 }-y_{i }$ is the lattice spacing at $y_i.$
The terms $\eta(r_i)$ and $\eta(r_i+r_j)$ above and in the following
should be understood to be zero
for $i\notin \{-M,\dots,M\}$ or $j\notin \{-M,\dots,M\}.$

We now assume that the atoms are also subject to an external force,
 $\tilde f_{i}(y_i),$
that is obtained from an external potential energy of the form
\begin{equation*}
\Po^a(\yb)=\sum_{i=-M}^{M+1}\Po^a_i(y_i),
\end{equation*}
so
\begin{equation}
\label{totalp}
\tilde f_{i}(y_i)=-\frac{\partial \Po^a(\yb)}{\partial y_{i }}
=-\frac{\partial \Po^a_i(y_i)}{\partial y_{i }}.
\end{equation}
For example, such external forces may model the interaction of
the one-dimensional
chain with atoms in layers above and below the chain, as in the Frenkel
Kontorova model~\cite{marder}.

We then have the equilibrium equations
$$
F_{i}^a(\yb)+\tilde f_{i}(y_i)=0,\qquad i=-M,\dots,M+1.
$$

\subsection{The Constrained Atomistic Quasicontinuum Approximation}
One can reduce the degrees of freedom in the atomistic model by linearly
interpolating the positions of the atoms between a set of
representative atoms (see
\cite{ortnersuli} for an analysis in this case).
We introduce representative
atoms with positions $z_{j}$ such that
\begin{equation} \label{lat}
z_{j}=y_{\ell_j}\quad\text{for }j=-N,\dots,N+1,
\end{equation}
where $\ell_{-N}=-M,$ $\ell_{N+1}=M+1,$ and $\ell_j<\ell_{j+1}.$  We then let
$\nu_j = \ell_{j+1} - \ell_{j}$ denote the number of atoms between $z_{j}$
and $z_{j+1}$ (see Figure \ref{contchain}). We now have that
\begin{equation*}
r_j=\frac{z_{j+1} - z_{j}}{\nu_j}
\end{equation*}
is the
distance separating $\nu_j$ equally spaced atoms between
$z_j$ and $z_{j+1},$ and
we have the conservation of mass equation
\begin{equation} \label{conserve}
\sum_{j=-N}^{N} \nu_j=\sum_{j=-N}^{N}\left( \ell_{j+1} - \ell_{j}\right)=
\ell_{N+1}-\ell_{-N}=2 M + 1.
\end{equation}
\begin{figure}[htb]
\includegraphics{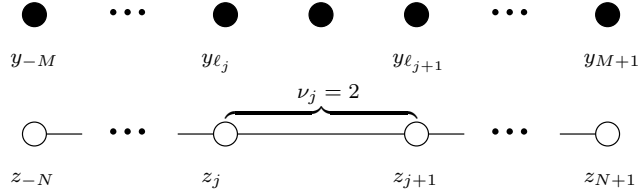}
\caption{\label{contchain} A coarsening of the atomistic chain.}
\end{figure}

In the constrained atomistic model,
the positions of atoms between $z_{j}$ and $z_{j+1}$
are linearly interpolated as
\begin{equation}\label{conat}
y_{\ell_{j} + i } =
y_{\ell_{j+1}-(\nu_j- i) } =
\frac{\nu_j - i}{\nu_j} z_{j} +
\frac{i}{\nu_j} z_{j+1} \quad\text{for }0\leq i \leq \nu_j,
\end{equation}
and we define the total internal energy in terms of
$\zb = (z_{-N},\dots,z_{N+1}) \in \Real^{2N+2}$ to
again be the interaction enegy of all atoms in the chain,
computed according to \eqref{atom2}, giving
\begin{equation}
\label{con}
\E^c(\zb) =  \sum_{j=-M}^{M+1} \E_j^a(\yb(\zb)).
\end{equation}

So far, we have reduced the degrees of freedom necessary for denoting
the atomistic positions, but the total energy is still computed as a sum of
energy contributions from all of the atomistic degrees of freedom.
However, all nearest-neighbor and next-nearest-neighbor contributions
for atoms interpolated between a pair of represenative atoms are
identical due to the uniform spacing, so if $y_i$ does not denote
a representative atom we have that
\begin{equation*}
\phi(y_{i+1}-y_{i-1}) = \phi(2 (y_{i+1} - y_i)) = \phi(2(y_{i}-y_{i-1}))
= \phi(2 r_j) \quad \text{ for }\ell_j < i < \ell_{j+1}.
\end{equation*}
These interactions account for all of the energy contributions except for
next-nearest-neighbor interactions that straddle a representative atom.
We can treat these interactions by observing that
\begin{equation*}
\begin{split}
\phi(y_{i +1}&-y_{i-1}) -\half \phi(2 (y_{i}-y_{i-1}))
-\half \phi(2 (y_{i+1} - y_{i})) \\ &=
\begin{cases}
\phi(r_{j-1}+r_j)-\half \phi(2r_{j-1})- \half \phi(2r_j) ,& \text{ if }y_i=z_j\text{ for some }j=-N,\dots,N+1,\\
0,&\text{ if }y_i\ne z_j\text{ for all }j=-N,\dots,N+1.
\end{cases}
\end{split}
\end{equation*}

We can therefore partition the total energy into the energy
of the region between
$z_{j}$ and $z_{j+1}$ for each $j=-N,\dots,N,$ plus interfacial energy terms
that account for interactions that straddle a representative atom.  We have
from \eqref{con} that
\begin{equation}
\begin{split}
\label{constrained_energy}
\E^c(\zb) &= \sum_{i=-M}^{M} \phi(y_{i+1}-y_{i})
+\negthickspace\sum_{i=-M+1}^{M} \phi(y_{i+1}-y_{i-1})\\
&= \sum_{i=-M}^{M} \left[\phi(y_{i+1}-y_{i}) + \phi(2(y_{i+1} - y_{i}))\right]\\
&\qquad\qquad +\negthickspace\sum_{i=-M+1}^{M} \left[\phi(y_{i+1}-y_{i-1})
-\half \phi(2 (y_{i}-y_{i-1})) -\half \phi(2 (y_{i+1} - y_{i}))\right]\\
&\qquad\qquad -\half\phi(2(y_{-M+1}-y_{-M}))-\half\phi(2(y_{M+1}-y_{M}))\\
&= \sum_{j=-N}^{N}
  \nu_j \hat\phi\left(r_j\right)
  +\sum_{j=-N}^{N+1}\Su_j\left(r_{j-1},
  r_j\right),
\end{split}
\end{equation}
where
\begin{equation}
\label{energydensity}
\hat\phi(r)=\phi(r)+\phi(2r)\quad\text{for }r>0
\end{equation}
and where
\begin{equation*}
\begin{split}
\Su_{-N}(r_{-N})=&-\half\phi(2r_{-N}),\\
\Su_j(r_{j-1},r_j) =& - \half\phi(2r_{j-1}) +  \phi(r_{j-1} + r_j)
    - \half\phi(2r_j),\qquad j=-N+1,\dots,N,\\
\Su_{N+1}(r_N)=&-\half\phi(2r_N).
\end{split}
\end{equation*}
We note that $\hat\phi(r)=\phi(r)+\phi(2r)$
is the energy per atom for an
infinite atomistic chain with the uniform lattice spacing $y_{i+1}-y_{i}=r$ for all
$-\infty<i<\infty,$ and that $\Su_j(r_{j-1},r_j)$
can be considered to be a surface energy at $z_{-N}$ and $z_{N+1}$
and to be an interfacial energy at $z_j$ for $j=-N+1,\dots,N.$
We observe that $\Su_j(r_{j-1},r_j)$ is a second divided difference for
$\phi(r)$ about $r=r_{j-1} + r_j$ with increment $r_j-r_{j-1},$ so
\[
\Su_j(r_{j-1},r_j) = - \half\phi''(r_{j-1} + r_j)(r_j-r_{j-1})^2
    +O(|r_j-r_{j-1}|^4),\qquad j=-N+1,\dots,N.
    \]
We will see that if $\phi(r)$ satisfies the
assumptions given in Section~\ref{exist}, then the convexity condition
\[
\Su_j(r_{j-1},r_j)>0,\qquad j=-N+1,\dots,N,
\]
holds for the range of $\rb=(r_{-N},\dots,r_N)\in\Real^{2N+1}$
defined by \eqref{what} in
Theorem~\ref{thm:wp1} and Corollary ~\ref{cor:wp1}
where solutions of the force-based quasicontinuum equilibrium
equations are shown to exist and where the iterates
of the ghost force iteration reside.

\subsection{The Local Quasicontinuum Approximation}

If we neglect the surface and interfacial energy terms, $\Su_j,$ in
\eqref{constrained_energy}, then we obtain the
local quasicontinuum approximation~\cite{tadmor_miller_qc_overview}
\[
\E^L(\zb) = \sum_{j=-N}^{N+1} \E_{j}^L(\zb),
\]
where
\begin{equation}\label{sss}
\E_{j}^L(\zb) = \half \left[\hat \phi(r_{j}) \nu_j +\hat \phi(r_{j-1})
    \nu_{j-1} \right].
\end{equation}
To treat the boundary terms consistently, we set $\nu_{-N-1}=\nu_{N+1}=1$
and $\hat\phi(r_{-N-1})=\hat\phi(r_{N+1})=0.$

We remark that the representative atoms need not be placed at atomistic
sites as described by \eqref{lat}.  The local quasicontinuum approximation only requires
that $z_{j}<z_{j+1}$ for $j=-N,\dots,N$ and that
the $\nu_j$ are positive and satisfy the conservation of mass
condition \eqref{conserve}.
The approximation
can be generalized to higher space dimensions by using the
Cauchy-Born rule~\cite{tadmor_miller_qc_overview}.
The local quasicontinuum approximation is computationally simpler than the
constrained atomistic quasicontinuum approximation,
especially in higher dimensions where the computation of the
interfacial energy becomes
expensive.  In the following, we will sometimes refer to the
local quasicontinuum approximation as the continuum approximation.

The force on a
representative atom at $z_{j}$ for
$j=-N+1,\dots,N$ is given for the local quasicontinuum approximation by
\begin{equation}\label{continuumforce}
\begin{split}
F_{j}^L(\zb) &=
-\frac{\partial}{\partial z_{j }}
\left[\hat\phi\left(\frac{z_{j + 1}-z_{j }}{\nu_{j}}\right)  \nu_{j}
+ \hat\phi\left(\frac{z_{j }-z_{j-1 }}{\nu_{j-1}}\right)\nu_{j-1}\right]
\\
&=-\frac{\partial}{\partial z_{j }}
\left[ \phi\left(\frac{z_{j + 1}-z_{j }}{ \nu_{j}}
\right)  \nu_{j}
+ \phi\left(\frac{2(z_{j + 1}-z_{j })}{ \nu_{j}}\right)
  \nu_{j} \right. \\
&\qquad\qquad+ \left. \phi\left(\frac{z_{j }-z_{j-1 }}{\nu_{j-1}}\right)
 \nu_{j-1}
+ \phi\left(\frac{2(z_{j }-z_{j-1})}{ \nu_{j-1}}\right) \nu_{j-1}
\right]  \\
&= \left[\eta(r_{j}) + 2 \eta(2r_{j})\right]
- \left[\eta(r_{j-1}) + 2 \eta(2 r_{j-1})\right],
\end{split}
\end{equation}
where again
\[
r_j = \frac{(z_{j+1}-z_{j })}{ \nu_j}
\]
is the lattice constant for the atoms between $z_j$ and $z_{j+1}.$
  In the above,  we see the ``local''
nature of the approximation, as the force on a degree of freedom is
determined only by the positions of adjacent degrees of freedom and
no long-rage interactions occur.

We can similarly compute the force on the boundary atoms, noting the one-sided
nature of $\E_{j}^L(\zb)$ for the boundary atoms,  by
\begin{equation}\label{continuumforceboundaryright}
\begin{split}
F_{-N}^L(\zb) &=
-\frac{\partial}{\partial z_{-N }}
\left[\hat\phi\left(\frac{z_{-N+1}-z_{-N }}{\nu_{-N}}\right) \nu_{-N}\right] \\
&=-\frac{\partial}{\partial z_{-N }}
\left[ \phi\left(\frac{z_{-N+1}-z_{-N}}{\nu_{-N}}\right)
+ \phi\left(\frac{2(z_{-N+1}-z_{-N})}{\nu_{-N}}\right)
\right] \nu_{-N} \\
&= \left[\eta(r_{-N}) + 2 \eta(2 r_{-N})\right],
\end{split}
\end{equation}
\begin{equation}\label{continuumforceboundaryleft}
\begin{split}
F_{N+1}^L(\zb) &=
-\frac{\partial}{\partial z_{N+1}}
\left[\hat\phi\left(\frac{z_{N+1}-z_{N }}{\nu_N} \right)\nu_N\right] \\
&=-\frac{\partial}{\partial z_{N+1}}
\left[ \phi\left(\frac{z_{N+1}-z_{N }}{ \nu_N}\right)
+ \phi\left(\frac{2(z_{N+1}-z_{N })}{\nu_N}\right)
\right] \nu_N \\
&= - \left[\eta(r_N) + 2 \eta(2 r_N)\right].
\end{split}
\end{equation}

We note that the local quasicontinuum energy, $\E^L(\zb),$
and the forces, $F_{j}^L(\zb),$ depend only on
$\rb=(r_{-N},\dots,r_N),$ and we will denote the
dependance by $\E^L(\rb)$ and $F_{j}^L(\rb)$
without introducing distinct functions.

We can also derive a local quasicontinuum approximation for the
external potential, $\Po^L(\zb),$ by
setting
\begin{equation}\label{pot}
\Po^L(\zb) = \Po^a(\yb(\zb)).
\end{equation}
By \eqref{totalp} and \eqref{conat}, the external force on the
representative atom at position $z_j$ is
\begin{equation*}
\begin{split}
f_{j}(\zb)=-\frac{\partial \Po^L (\zb)}{\partial z_{j }}
&=-\sum_{i=0}^{\nu_{j-1}}
  \left(\frac{\nu_{j-1} -i}{\nu_{j-1}}\right)
   \frac{\partial \Po^a(\yb(\zb))}{\partial y_{\ell_j-i}}
  -\sum_{i=1}^{\nu_{j}} \left(\frac{\nu_{j} - i}{\nu_{j}}\right)
  \frac{\partial \Po^a (\yb(\zb))}{\partial y_{\ell_j+i}}\\
&=\sum_{i=0}^{\nu_{j-1}}
  \left(\frac{\nu_{j-1} -i}{\nu_{j-1}}\right)
  \tilde f_{\ell_j-i}\left(y_{\ell_j-i}(\zb)\right)
  + \sum_{i=1}^{\nu_{j}} \left(\frac{\nu_{j} - i}{\nu_{j}}\right)
  \tilde f_{\ell_j+i}\left(y_{\ell_j+i}(\zb)\right).
\end{split}
\end{equation*}
It follows from the linear interpolation~\eqref{conat} that
\[
f_{j}(\zb)=f_{j}(z_{j-1},z_j,z_{j+1}).
\]

We shall assume in our analysis that the external forces,
$\tilde f_i,$ are independent of $\yb.$  In this case,
the local quasicontinuum forces, $f_j,$ are independent of
$\zb$ and
\begin{equation}\label{quasif}
f_{j}=\sum_{i=0}^{\nu_{j-1}}
  \left(\frac{\nu_{j-1} -i}{\nu_{j-1}}\right)\tilde f_{\ell_j-i}
  +\sum_{i=1}^{\nu_{j}} \left(\frac{\nu_{j} - i}{\nu_{j}}\right)
\tilde f_{\ell_j+i},\qquad j=-N,\dots,N+1.
\end{equation}
We consider any term $\tilde f_j$ to be
zero if $j \notin \{-M,\dots,M+1\}$ and any term $\nu_j$ to be one if
$j \notin \{-N,\dots,N\}.$

\subsection{The Energy-Based Quasicontinuum Approximation}

To describe the energy-based quasicontinuum approximation
~\cite{tadmor_miller_qc_overview}, we again introduce representative
atoms with positions $z_{j}$ for $j=-N,\dots,N+1,$ where $z_{j}<z_{j+1}.$
Each representative atom is considered to be an ``atomistic'' or
``continuum'' degree of freedom and will contribute either
$\E_j^a(\zb)$ or $\E_j^L(\zb)$ to the total
internal energy according to the atomistic model~\eqref{atompart}
or the local
quasicontinuum approximation~\eqref{sss}, respectively.

In applications, atomistic degrees of freedom are used in regions of
interest where highly non-uniform behavior is expected.  Continuum
regions surround this, gradually coarsening by increasing $\nu_j$
in regions with slowly varying strain.  For simplicity of exposition,
we will consider an approximation with a single atomistic region,
symmetrically surrounded by continuum regions large enough so
that no atomistic degrees of freedom interact with the surface atoms
through nearest-neighbor or next-nearest-neighbor interactions.

We denote the representative atom positions by $z_j$ and define the range
$j = -K+1,\dots,K$ to be atomistic sites and the ranges $j=-N,\dots, -K$ and
$K+1,\dots, N+1$ to be continuum sites.  Therefore, the total quasicontinuum
energy, $\E^{QC}(\rb),$ for the chain is given by
\begin{equation}\label{total}
\E^{QC}(\rb) = \sum_{j = -N}^{-K} \E_{j }^L(\rb) +
\sum_{j=-K+1}^{K} \E_{j}^a(\rb) +
\sum_{j=K+1}^{N+1} \E_{j }^L(\rb),
\end{equation}
where $\E_{j}^a(\rb)$ is defined in \eqref{atompart} and $\E_{j }^L(\rb)$ is
defined in \eqref{sss}.
We assume that $\nu_j=1$ for $j=-K-1,\dots, K+1.$ This guarantees that
$\nu_j=1$ within the next-nearest-neighbor
cutoff radius of any atomistic site and enables a
seamless transition to the continuum approximation.

\begin{figure}[htb]
\includegraphics{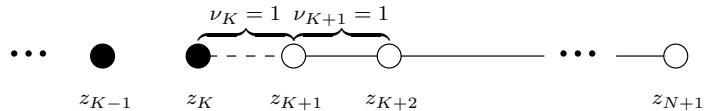}
\caption{One end of the quasicontinuum chain, highlighting the interface.
Filled circles are atomistic representative atoms, whereas the
unfilled circles are
continuum representative atoms.}
\end{figure}

The force, $F_{j}^{QCE}(\rb),$ at $z_j$ is then given by
\begin{equation}\label{colu}
\begin{split}
F_{j}^{QCE}(\rb) &=-\frac{\partial \E^{QC}}{\partial z_j}(\rb)
=-\sum_{\ell=-N}^N \frac{\partial \E^{QC}}{\partial r_\ell}(\rb)
\frac{\partial r_\ell}{\partial z_j}\\
&=\begin{cases}
F_j^L(\rb), &  -N \leq j \leq -K-2, \\
F_{-K-1}^L(\rb) + \half \eta(r_{-K-1} + r_{-K}), & j = -K-1, \\
F_{-K}^L(\rb) - \eta(2 r_{-K}) + \half \eta(r_{-K} + r_{-K+1}), & j = -K,\\
F_{-K+1}^a(\rb) - \eta(2 r_{-K}) + \half \eta(r_{-K-1} + r_{-K}), & j = -K+1,\\
F_{-K+2}^a(\rb) + \half \eta(r_{-K} + r_{-K+1}), & j = -K + 2, \\
F_j^a(\rb), &  -K+3 \leq j \leq K-2,\\
F_{K-1}^a(\rb) - \half \eta(r_{K-1} + r_{K}), & j = K-1,\\
F_{K}^a(\rb)   + \eta(2 r_K) - \half \eta(r_K + r_{K+1}), & j = K, \\
F_{K+1}^L(\rb) + \eta(2 r_K) - \half \eta(r_{K-1} + r_K), & j = K+1, \\
F_{K+2}^L(\rb) - \half \eta(r_{K}+r_{K+1}), & j = K+2, \\
F_j^L(\rb), & K+3 \leq j \leq N+1.
\end{cases}
\end{split}
\end{equation}

In the above expression, we notice that in the large ranges interior to
the atomistic and continuum regions the forces are exactly those from the
individual models,
namely, either $F_j^a(\rb)$ or $F_j^L(\rb).$
Near the atomistic to continuum interface, there are additions
to these force terms which contain non-physical
``ghost'' forces.

\begin{figure}[htb]
\includegraphics{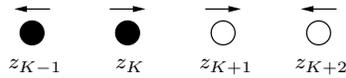}
\caption{\label{resultant}Direction of nonzero forces for a
uniform configuration.}
\end{figure}
To see this, we consider the forces on the representative
atoms when the lattice spacings are uniform, that is, when
$r_j=a$ for $j=-N,\dots,N,$ or equivalently when $\rb=\az=(a,\dots,a)\in\Real^{2N+1}.$
We observe that the forces computed
according to the atomistic model \eqref{atomforce} or the continuum model
\eqref{continuumforce} are zero
except at the ends of the chain, that is, $F_j^L(\az)=0$ for
$j=-N+1,\dots,N$ and $F_j^a(\az)=0$ for
$j=-N+2,\dots,N-1.$  We then have from \eqref{colu} that
\begin{equation*}
F_{j}^{QCE}(\az)=
\begin{cases}
\left[\eta(a) + 2 \eta(2a)\right], & j = -N,\\
0, &  -N+1 \leq j \leq -K-2, \\
  \half \eta(2a), & j = -K-1, \\
- \half \eta(2a), & j = -K,\\
- \half \eta(2a), & j = -K+1,\\
  \half \eta(2a), & j = -K + 2, \\
0, &  -K+3 \leq j \leq K-2,\\
- \half \eta(2a), & j = K-1,\\
  \half \eta(2a), & j = K, \\
  \half \eta(2a), & j = K+1, \\
- \half \eta(2a), & j = K+2, \\
0, & K+3 \leq j \leq N, \\
-\left[\eta(a)+2 \eta(2a)\right], & j = N+1.
\end{cases}
\end{equation*}
However, a continuum chain with
uniform lattice spacings has forces only at the surfaces.
Figure \ref{resultant} shows the ghost forces at one of the
interfaces.  Figure \ref{forces} shows the origin of out of
balance forces for a single representative atom.  We note that while the first
neighbor terms are balanced, the second neighbor terms are not
balanced due to the fact that the continuum site does not contribute
any second-neighbor interactions.
\begin{figure}[htb]
\includegraphics{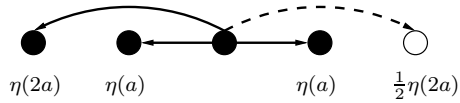}
\caption{\label{forces}Imbalance of forces on an atomistic representative atom
near the interface.}
\end{figure}

The force-based quasicontinuum approximation described in the
next subsection corrects these ghost forces.

\subsection{The Force-Based Quasicontinuum Approximation}
The force-based quasicontinnum approximation~\cite{tadmor_miller_qc_overview}
corrects the nonphysical forces described in the previous subsection.
The forces on the representative
atoms in the interior of the atomistic and continuum regions are
defined by \eqref{atomforce} and \eqref{continuumforce}, but
the forces on the representative atoms near the atomistic-continuum
transition given by the energy-based quasicontinuum method must be
modified to remove the non-physical terms.

In the force-based quasicontinuum approximation, we again partition into
atomistic and continuum representative atoms, where the force on a
representative atom is
the force that would result on it if the approximation was entirely of its
respective
type (atomistic or continuum). With this convention, a continuum
representative atom only interacts with adjacent degrees of freedom
regardless of how close any atomistic sites may be.  We will see
that the tradeoff for this simple philosophy is that the forces
are not conservative, that is, they cannot be derived from an energy.

We now model the forces on the representative atoms for
$j=-K+1,\dots, K$ where $K<N-1$ by the atomistic model
\eqref{atomforce} and the forces on the representative
atoms for $j=-N,\dots, -K$ and $j=K+1,\dots, N+1$ by
the local quasicontinuum approximation
~\eqref{continuumforce}-\eqref{continuumforceboundaryleft}.
The forces on the representative atoms
are then given by
\begin{align*}
F_{j}^{QCF}(\rb) &=
\begin{cases}
F_j^L(\rb), & -N \leq j \leq -K,\\
F_j^a(\rb),&-K+1 \leq j \leq K,\\
F_j^L(\rb),&K+1 \leq j \leq N+1,
\end{cases}\\
&=
\begin{cases}
\left[\eta(r_{-N} ) + 2 \eta(2 r_{-N} ) \right], &j = -N,\\
\left[\eta(r_{j}) + 2 \eta(2 r_{j})\right] -
  \left[\eta(r_{j-1}) + 2 \eta(2 r_{j-1})\right], & -N+1 \leq j \leq -K,\\
\left[\eta(r_{j}) + \eta(r_{j} + r_{j+1})\right] -
  \left[\eta(r_{j-1}) + \eta(r_{j-1} + r_{j-2})\right],&-K+1 \leq j \leq K,\\
\left[\eta(r_{j}) + 2 \eta(2 r_{j})\right] -
  \left[\eta(r_{j-1}) + 2 \eta(2 r_{j-1})\right],&K+1 \leq j \leq N,\\
-\left[\eta(r_N ) + 2 \eta(2 r_N )\right], & j = N+1.
\end{cases}
\end{align*}

The force-based quasicontinuum formulation has the desired property
that if we take a uniform configuration, $\az=(a,\dots,a)\in\Real^{2N+1}$,
we have that
\[
F_{j}^{QCF}(\az)= 0, \qquad j = -N+1,\dots,N,
\]
and on the boundary we get values of equal magnitude and opposite
signs.  Thus, the equilibrium equations have a uniform solution, $\az$,
whenever the external forces are applied at the boundary, that is, a chain in
uniform tension or compression.

However, the solution to the equilibrium equations for the
force-based quasicontinuum method,
\begin{equation}
\label{equilq}
F_{j}^{QCF}(\rb)+f_{j}=0,\qquad j=-N,\dots,N+1,
\end{equation}
cannot be obtained from the minimization of an energy since
$F_{j}^{QCF}(\rb)$ is a nonconservative force.  To see this,
we observe that the
forces given by the force-based quasicontinuum approximation,
$F_{j}^{QCF}(\rb),$ are not the differential of an energy, $\E^{QCF}(\rb),$
since $F_{j}^{QCF}(\rb)$ is not a closed form~\cite{fleming}.
We can see this by noting that
\begin{equation*}
\begin{split}
\frac{\partial F_{K}^{QCF}}{\partial z_{K+1}}(\rb)&=\eta'(r_K),\\
\frac{\partial F_{K+1}^{QCF}}{\partial z_{K}}(\rb)&=\eta'(r_K)+4\eta'(2 r_K),
\end{split}
\end{equation*}
implies that
\[
\frac{\partial F_{K}^{QCF}}{\partial z_{K+1}}(\rb)\ne
\frac{\partial F_{K+1}^{QCF}}{\partial z_{K}}(\rb).
\]

\subsection{The Ghost Force Iteration}

In the quasicontinuum method of \cite{tadmor_miller_qc_overview},
the forces, $F^{QCF}(\rb), $ are split into the force from the
energy-based quasicontinuum approximation, $F^{QCE}(\rb),$ and
ghost force corrections defined by
\begin{equation} \label{gfdef}
F^G_j(\rb) =F^{QCF}_j(\rb) - F^{QCE}_j(\rb)=
\begin{cases}
0, &  -N \leq j \leq -K-2, \\
- \half \eta(r_{-K-1} + r_{-K}), & j = -K-1, \\
+ \eta(2 r_{-K}) - \half \eta(r_{-K} + r_{-K+1}), & j = -K,\\
+ \eta(2 r_{-K}) - \half \eta(r_{-K-1} + r_{-K}), & j = -K+1,\\
- \half \eta(r_{-K} + r_{-K+1}), & j = -K + 2, \\
0, &  -K+3 \leq j \leq K-2,\\
+ \half \eta(r_{K-1} + r_{K}), & j = K-1,\\
- \eta(2 r_K) + \half \eta(r_K + r_{K+1}), & j = K, \\
- \eta(2 r_K) + \half \eta(r_{K-1} + r_K), & j = K+1, \\
+ \half \eta(r_{K}+r_{K+1}), & j = K+2, \\
0, & K+3 \leq j \leq N+1.
\end{cases}
\end{equation}
The forces $ F^G_j(\rb)$ act as a model correction near the atomistic-continuum interfaces
to enforce the convention that each representative atom
 has forces acting on it
as though it were surrounded by representative atoms of the same type.

In \cite{tadmor_miller_qc_overview},
the solution to the equilibrium equations
\eqref{equilq} is obtained by solving the iteration
\begin{equation}
\label{iter1}
F_j^{QCE}(\rb^{n+1}) + F_j^G(\rb^{n}) + f_j = 0, \qquad j = -N,\dots,N+1,
\end{equation}
by using a conjugate gradient method to compute $\rb^{n+1}$ from $\rb^{n}.$
If the sequence of solutions $\{{\rb}^n\}$
converges, then the iterative limit $\rb$ satisfies
the equilibrium equations~\eqref{equilq}.

We have split the internal forces of the force-based quasicontinuum approximation,
$F^{QCF}_j(\rb),$ into a conservative force,
$F^{QCE}_j(\rb),$ and correction, $F^{G}_j(\rb),$ to define the iterative method
above \eqref{iter1}.
The conservative force, $F^{QCE}_j(\rb),$ was defined by the quasicontinuum energy,
$\E^{QC}(\rb),$ given by \eqref{total}, and the correction, $F^{G}_j(\rb),$
was defined simply by $F^{G2}_j(\rb) =F^{QCF}_j(\rb) - F^{QCE2}_j(\rb).$
 However, we are free to try to improve the rate of
convergence of $\rb^n\to\rb$
by constructing an energy $\E^{QCE2}(\rb)$ different from $\E^{QCE}(\rb)$ and
then solving the iteration
\begin{equation*}
F_j^{QCE2}(\rb^{n+1}) + F_j^{G2}(\rb^{n}) + f_j = 0, \qquad j = -N,\dots,N+1,
\end{equation*}
where $F_j^{QCE2}(\rb)$ is now the force that is derived from $\E^{QCE2}(\rb)$ and
$F^{G2}_j(\rb) =F^{QCF}_j(\rb) - F^{QCE2}_j(\rb).$

\section{Conjugate Forces}
\label{integ}

In this section, we simplify the analysis by reformulating the
force-based quasicontinuum equilibrium equations
\eqref{equilq} in terms of forces
conjugate to the distance between representative atoms,
\[
R_j=z_{j+1}-z_j=\nu_j{r_j},\qquad j=-N,\dots,N,
\]
 rather than conjugate to the positions of the representative atoms, $z_j.$
We will use the notation $\Rb=(R_{-N},\dots,R_N)\in\Real^{2N+1}.$
Through this technique, we will be able to derive equations that are decoupled
inside the continuum regions and are the sum of tridiagonal terms
and nonlocal interfacial terms in the atomistic region.

\subsection{The Internal Conjugate Force}

We define the internal conjugate force for
the energy-based quasicontinuum approximation by
\begin{equation}\label{not}
\psi^E_j(\rb)=\frac{\partial \E^{QC}}{\partial R_j}\left(\rb(\Rb)\right)
=\frac1{\nu_j}\frac{\partial \E^{QC}}{\partial r_j}(\rb), \qquad
j=-N,\dots,N.
\end{equation}
We note that we have found it convenient to define
$\psi^E_j(\rb)$ as the negative of the usual convention
for a conjugate force.
We next derive the following relation between $F^{QCE}_j(\rb)$
and the internal conjugate force $\psi^E_j(\rb):$
\begin{equation}\label{for}
\begin{split}
F^{QCE}_j(\rb)&=-\frac{\partial \E^{QC}}{\partial z_j}(\rb)
=-\frac{\partial \E^{QC}}{\partial R_j}\left(\rb(\Rb)\right)
\frac{\partial R_j}{\partial z_j}
-\frac{\partial \E^{QC}}{\partial R_{j-1}}\left(\rb(\Rb)\right)\frac{\partial R_{j-1}}{\partial z_j}
\\
&=\psi^E_{j}(\rb) - \psi^E_{j-1}(\rb),\qquad j=-N,\dots,N+1,
\end{split}
\end{equation}
where we set
\[
\psi^E_{-N-1}(\rb)=\psi^E_{N+1}(\rb)=0.
\]
We can sum the forces from the left
of the chain and use the preceding equation~\eqref{for}
to obtain that
\begin{equation}\label{sumstressz}
\psi^E_j(\rb) = \sum_{i=-N}^j F^{QCE}_{i}(\rb),\qquad j=-N,\dots,N+1.
\end{equation}
We can derive from either \eqref{not} or \eqref{sumstressz} that
\begin{equation}
\label{psiE}
\psi^E_j(\rb) =
\begin{cases}
0,&j=-N-1\\
  \eta(r_j) + 2 \eta(2r_j), &  -N \leq j \leq -K-2,\\
  \eta(r_j) + 2 \eta(2r_j) + \half \eta(r_j +r_{j+1}), & j = -K-1, \\
  \eta(r_j) + \half \eta(r_j+r_{j-1}) + \half \eta(r_j+r_{j+1}) +
      \eta(2r_j), & j = -K,\\
  \eta(r_j) + \half \eta(r_j+r_{j-1}) + \eta(r_j+r_{j+1}),  & j = -K+1, \\
  \eta(r_j) + \eta(r_j+r_{j-1}) + \eta(r_j+r_{j+1}),
       & -K + 2 \leq j \leq K - 2, \\
  \eta(r_j) + \eta(r_j+r_{j-1}) + \half \eta(r_j+r_{j+1}),  & j = K-1, \\
  \eta(r_j) + \half \eta(r_j+r_{j-1}) + \half \eta(r_j+r_{j+1}) +
      \eta(2r_j), & j = K,\\
  \eta(r_j) + 2 \eta(2r_j) + \half \eta(r_j +r_{j-1}), & j = K+1, \\
  \eta(r_j) + 2 \eta(2r_j), &  K+2 \leq j \leq N,\\
  0,&j=N+1.
\end{cases}
\end{equation}

We cannot properly derive an internal conjugate force for the force-based
quasicontinuum approximation since it is not a conservative force.
However, we will find it convenient to
define an internal conjugate force for
the force-based
quasicontinuum approximation by following \eqref{sumstressz} setting
\[
\psi^F_{-N-1}(\rb)=0
\]
and
\begin{equation}\label{sumstress}
\psi^F_j(\rb) = \sum_{i=-N}^j F^{QCF}_{i}(\rb),\qquad j=-N,\dots,N+1.
\end{equation}
We can then obtain that
\begin{equation}\label{diffstress}
F^{QCF}_{j}(\rb) = \psi^F_{j}(\rb) - \psi^F_{j-1}(\rb),\qquad j=-N,\dots,N+1.
\end{equation}

 We
have the following closed form expressions for the internal conjugate force
$\psi^F_j(\rb):$
\begin{equation}
\label{stress1}
\psi^F_j(\rb) =
\begin{cases}
0,& j=-N-1,\\
\eta(r_j) + 2 \eta(2 r_j), & -N \leq j \leq -K ,\\
\eta(r_j) + \eta(r_j + r_{j-1}) + \eta(r_j + r_{j+1}) & \\
\qquad + [2\eta(2r_{-K}) - \eta(r_{-K} + r_{-K-1}) - \eta(r_{-K} + r_{-K+1})],
    & -K+1 \leq j \leq K ,\\
\eta(r_j) + 2 \eta(2 r_j) & \\
\qquad + [2\eta(2r_{-K}) - \eta(r_{-K} + r_{-K-1}) - \eta(r_{-K} + r_{-K+1})]
    & \\
\qquad- [2\eta(2r_{K}) - \eta(r_{K} + r_{K-1}) - \eta(r_{K} + r_{K+1})],
    & K+1 \leq j \leq N,\\
[2\eta(2r_{-K}) - \eta(r_{-K} + r_{-K-1}) - \eta(r_{-K} + r_{-K+1})]&\\
\qquad- [2\eta(2r_{K}) - \eta(r_{K} + r_{K-1}) - \eta(r_{K} + r_{K+1})],
    & j= N+1.
\end{cases}
\end{equation}

The internal conjugate force, $\psi^F_j(\rb),$ takes a simpler form when there is no
resultant quasicontinuum force, that is, when $\rb$ satisfies
\begin{equation}\label{resul}
\psi^F_{N+1}(\rb) = \sum_{j=-N}^{N+1} F_{j}^{QCF}(\rb) = 0,
\end{equation}
which is equivalent to
\begin{equation}
\label{noresultant}
\begin{split}
2\eta(2r_{-K}) &- \eta(r_{-K} + r_{-K-1}) - \eta(r_{-K} + r_{-K+1}) \\
   &= 2\eta(2r_{K}) - \eta(r_{K} + r_{K-1}) - \eta(r_{K} + r_{K+1}).
\end{split}
\end{equation}
We note that \eqref{noresultant} is satisfied when
$\rb$ is symmetric, that is, when $r_{-j} = r_j$ for $j=1,\dots,N.$
Let us now define $\hat\psi_j^F(\rb)$ on all
of $\Real^{2N+1}$ as a symmetric extension of
$\psi^F_j(\rb),$ with equality whenever $\rb$ satisfies~\eqref{noresultant}.
This leads to $\hat\psi^F_j(\rb)$ having a more symmetric form,
\begin{align}
\label{stress2}
\hat\psi^F_j(\rb) =
\begin{cases}
0,&j=-N-1,\\
\eta(r_j) + 2 \eta(2 r_j),
       & -N \leq j \leq -K ,\\
\eta(r_j) + \eta(r_j+r_{j-1}) + \eta(r_j+r_{j+1}) &\\
\qquad+ [2\eta(2r_{K}) - \eta(r_{K} + r_{K-1}) - \eta(r_{K} + r_{K+1})],
     & -K+1 \leq j \leq K-1,\\
\eta(r_j) + 2 \eta(2 r_j),
     & K\leq j \leq N,\\
0,& j=N+1.
\end{cases}
\end{align}
The intervals of definition changed slightly from \eqref{stress1},
which is one of the simplifications afforded by \eqref{noresultant}.
We use $\hat \psi^F_j(\rb)$ in our subsequent analysis, and in Section
\ref{equileq} we discuss the relation between using $\psi^F_j(\rb)$
and $\hat{\psi}^F_j(\rb)$ to solve the equilibrium equations~\eqref{equilq}.

We identify in \eqref{stress2} a continuum internal conjugate force
for $j=-N,\dots,-K$ and $j=K,\dots,N$ given by
\[
\eta(r_j) + 2 \eta(2 r_j)
\]
and an atomistic internal conjugate force for $j=-K+1,\dots,K-1$
given by
\[
\eta(r_j) + \eta(r_{j-1} + r_j) + \eta(r_j + r_{j+1}).
\]
We identify the remaining terms,
\[
2\eta(2r_{K}) - \eta(r_{K} + r_{K-1}) - \eta(r_{K} + r_{K+1}),
\]
for $j=-K+1,\dots,K-1$ as the nonlocal part of the internal conjugate force.

 For consistency, we wish to define $\hat{\psi}^E_j(\rb)$ as
we did for $\hat{\psi}^F_j(\rb).$  Since $\psi^E_j(\rb)$
is derived from the energy $\E^{QCE}(\rb),$ it has no resultant
force. Thus, we define $\hat{\psi}^E_j(\rb)=\psi^E_j(\rb).$

We can also derive a corresponding internal conjugate force, $\psi^G_j(\rb),$
by summing $F^G_j(\rb)$ as in~\eqref{sumstress}.
From the definition of $F^G_j(\rb)$ \eqref{gfdef}, we have that
\[
\psi^F_j(\rb)=\psi^E_j(\rb)+\psi^G_j(\rb),\qquad j=-N-1,\dots,N+1.
\]
We can define $\hat\psi^G_j(\rb)$ by
\[
\hat\psi^F_j(\rb)=\hat\psi^E_j(\rb)+\hat\psi^G_j(\rb),\qquad j=-N-1,\dots,N+1,
\]
and we can then check that
\begin{equation}
\label{psiG}
\hat \psi^G_j(\rb) =
\begin{cases}
0, & -N-1 \leq j \leq -K-2, \\
- \half \eta(r_{-K} + r_{-K-1}), & j = -K-1,\\
\eta(2r_{-K}) - \half \eta(r_{-K}+r_{-K+1})
  - \half \eta(r_{-K}+r_{-K-1}), & j = -K, \\
2\eta(2r_{-K}) - \half \eta(r_{-K}+r_{-K+1}) - \eta(r_{-K}+r_{-K-1}),
  & j = -K+1, \\
2\eta(2r_K) - \eta(r_K+r_{K-1}) - \eta(r_K+r_{K+1}),
  & -K+2 \leq j \leq K - 2, \\
2\eta(2r_K) - \half \eta(r_K+r_{K-1}) - \eta(r_K+r_{K+1}),   & j = K-1, \\
\eta(2r_K) - \half \eta(r_K+r_{K-1}) - \half \eta(r_K+r_{K+1}), & j = K, \\
- \half \eta(r_K + r_{K+1}), & j = K+1,\\
0, & K+2 \leq j \leq N+1.
\end{cases}
\end{equation}
If
$\rb$ satisfies the condition~\eqref{noresultant},
then we have that $\hat\psi^G_j(\rb)=\psi^G_j(\rb).$
We note that the nonlocal part of the
internal conjugate force is found in $\psi^G_j(\rb)$  and $\hat\psi^G_j(\rb).$

We finally observe that if $\rb$ is symmetric, then
$\hat\psi^F_j(\rb)=\psi^F_j(\rb),$
$\hat\psi^E_j(\rb)=\psi^E_j(\rb),$ and $\hat\psi^G_j(\rb)=\psi^G_j(\rb)$
are symmetric.

\subsection{The External Conjugate Force}
We recall that we are assuming in our analysis that the external forces,
$\tilde f_i,$ are independent of $\yb,$  and that consequently
the local quasicontinuum forces, $f_j,$ given in
\eqref{quasif} are independent of
$\zb.$  The external potential, $\Po^L(\zb),$ given in \eqref{pot}
thus has the form
\[
\Po^L(\zb)=-\sum_{j=-N}^{N+1}f_j z_j.
\]

We now assume that there is also no resultant force
from the external forces, so that
\begin{equation}\label{nores}
\sum_{j=-N}^{N+1} f_j = 0.
\end{equation}
It then follows that the external potential, $\Po^L(\zb),$ is also a function
of $\rb,$ and we can define the external conjugate force by
\begin{equation}\label{notz}
\Phi_j=-\frac{\partial \Po^L}{\partial R_j}\left(\rb(\Rb)\right)
=-\frac1{\nu_j}\frac{\partial \Po^L}{\partial r_j}(\rb), \qquad
j=-N,\dots,N.
\end{equation}

We next derive the following relation between $f_j$
and the external conjugate force $\Phi_j:$
\begin{equation}\label{forr}
\begin{split}
f_j&=-\frac{\partial \Po^L}{\partial z_j}(\rb)
=-\frac{\partial \Po^L}{\partial R_j}\left(\rb(\Rb)\right)\frac{\partial R_j}{\partial z_j}
-\frac{\partial \Po^L}{\partial R_{j-1}}\left(\rb(\Rb)\right)\frac{\partial R_{j-1}}{\partial z_j}
\\
&=-(\Phi_j - \Phi_{j-1}),\qquad j=-N,\dots,N+1,
\end{split}
\end{equation}
where we set
\[
\Phi_{-N-1}=\Phi_{N+1}=0.
\]
We can sum the external forces from the left
of the chain and use the preceding equation~\eqref{forr}
to obtain that
\begin{equation}\label{symm}
\Phi_j = -\sum_{i=-N}^j f_i,\qquad j=-N,\dots,N+1.
\end{equation}

If the
external forces, $f_{j},$ are anti-symmetric, that is,
\begin{equation}\label{forcesym}
f_{j+1}=-f_{-j},\qquad j=0,\dots,N,
\end{equation}
then
we can conclude from \eqref{symm} that the external conjugate force,
$\Phi_j,$ is symmetric about
$j=0,$ that is,
\begin{equation}\label{Phisym}
\Phi_{j}=\Phi_{-j},\qquad j=-N-1,\dots,N+1.
\end{equation}
We note that the external forces, $f_{j},$
are anti-symmetric
if the chain
is subject only to tensile or compressive loads of
equal magnitude, but opposite sign, at its ends.

\subsection{Equilibrium Equations}
\label{equileq}
It follows from  ~
\eqref{sumstress}, \eqref{diffstress}, \eqref{symm}, and \eqref{forr}
that the force-based quasicontinuum equilibrium equations
\begin{equation}
\label{equil}
F_{j}^{QCF}(\rb)+f_{j}=0,\qquad j=-N,\dots,N+1,
\end{equation}
are equivalent
to
\begin{equation*}
\psi^F_j(\rb) = \Phi_j,\qquad j=-N-1,\dots,N+1.
\end{equation*}
We recall that $\Phi_{N+1}=0$ since we have assumed that the
external forces satisfy the condition of no resultant force~\eqref{nores}.
Therefore, if $\rb$ satisfies the force-based quasicontinnum equilibrium equations
\eqref{equil}, then $\psi_{N+1}(\rb)=0,$ and by \eqref{resul} we have
$\psi^F_j(\rb)=\hat\psi_j^F(\rb)$ for
all $j=-N-1,\dots,N+1.$
Thus, we see that if $\rb$ satisfies
the force-based equilibrium equations~\eqref{equil}, then
\begin{equation*}
\hat\psi_j^F(\rb) = \Phi_j,\qquad j=-N-1,\dots,N+1.
\end{equation*}

We also recall that if $\rb$ is symmetric, then
\begin{equation*}
\psi^F_j(\rb)=\hat\psi_j^F(\rb),\qquad j=-N-1,\dots,N+1.
\end{equation*}
Hence, if $\rb$ is a symmetric solution of
\begin{equation*}
\hat\psi_j^F(\rb) = \Phi_j,\qquad j=-N-1,\dots,N+1,
\end{equation*}
then $\rb$ is a solution of
the force-based quasicontinuum equilibrium equations~\eqref{equil}.

If we sum the ghost force iteration equations
\begin{equation}
\label{iter}
F_j^{QCE}(\rb^{n+1}) + F_j^G(\rb^{n}) + f_j = 0, \qquad j = -N,\dots,N+1,
\end{equation}
then we also get
the following equivalent corresponding iterative method in terms
of the internal and external conjugate forces
\begin{equation*}
\psi^E_j(\rb^{n+1}) + \psi^G_j(\rb^{n}) = \Phi_j,\qquad j=-N-1,\dots,N+1.
\end{equation*}
We also have that if the sequence $\{ \rb^{n}\}$ satisfies the iteration equations
\eqref{iter}, then the sequence $\{\rb^{n}\}$ satisfies
\begin{equation}
\label{iterPsih}
\hat\psi^E_j(\rb^{n+1}) + \hat\psi^G_j(\rb^{n}) = \Phi_j,\qquad j=-N-1,\dots,N+1.
\end{equation}
We finally note that if the sequence $\{\rb^{n}\}$ is a
symmetric solution of \eqref{iterPsih},
then the sequence $\{\rb^{n}\}$ satisfies the ghost force iteration equations \eqref{iter}.

\section{Existence of Solutions to the Force-Based Quasicontinuum System}
\label{exist}

In this section, we will give conditions on
$\Pb=(\Phi_{-N},\dots,\Phi_N)\in\Real^{2N+1}$ such that
there exists a unique solution $\rb$ to
\begin{equation}\label{eb}
\hat{\mathbf\Psi}^F(\rb)=\Pb
\end{equation}
in a domain $\Omega\subset\Real^{2N+1},$ where
$\hat{\mathbf\Psi}^F(\rb)=\left(\hat{\psi}^F_{-N}(\rb),\dots,\hat{\psi}_N^F(\rb)\right)\in\Real^{2N+1}$
has the symmetric form given in \eqref{stress2}.
We note that we
ignore $\hat{\psi}^F_{-N-1},\,\hat{\psi}^F_{N+1},\,\Phi_{-N-1},\,\Phi_{N+1}$
in the above formulation since
$\hat{\psi}^F_{-N-1}=\hat{\psi}^F_{N+1}=\Phi_{-N-1}=\Phi_{N+1}=0$.
We conclude this section by showing that the unique
symmetric solution~$\rb$ to \eqref{eb}
for a symmetric $\Pb$ is a
solution of the force-based equilibrium equations \eqref{equil}.

For any uniform lattice spacing $\az = (a,\dots,a) \in \Real^{2N+1},$
we have a corresponding uniform
$\hat{\mathbf\Psi}^F(\az).$
  Under appropriate assumptions on the atomistic potential, $\phi(r),$
and the interatomic spacing, $\az,$ we will give an explicit
 neighborhood around $\az$
in which we have a unique solution to \eqref{eb} for
any external potential, $\Pb,$ in an explicit neighborhood of
$\hat{\mathbf\Psi}^F(\az).$
We will show that these assumptions on the atomistic potential, $\phi(r),$
are satisfied by the Lennard-Jones potential
\begin{equation}
\label{lj}
\phi(r) = \frac{1}{r^{12}} - \frac{2}{r^6}.
\end{equation}

\subsection{Assumptions on the Atomistic Potential,  $\phi(r)$}
We now give the assumptions on the atomistic potential, $\phi(r),$
that are required for our analysis.  We will assume that
$\phi(r) \in C^3\left((0,\infty\right)),$ and we recall that
$\eta(r)=\phi'(r)$ and $\hat\phi(r)=\phi(r)+\phi(2r).$
We define
\[
\hat\eta(r)=\hat\phi'(r),
\]
so
\[
\hat\eta(r)=\phi'(r)+2\phi'(2r)=\eta(r)+2\eta(2r).
\]

We assume that the atomistic potential, $\phi(r),$ satisfies the following
properties that are graphically displayed in
Figures \ref{phiplot} and \ref{hatphiplot}:
\begin{align}
\label{d1} &\eta'(r) > 0 \text{ for } 0 < r < \tilde r_1 \text{ and }
           \eta'(r) < 0 \text{ for } r > \tilde r_1,\\
\label{d2} &\eta''(r)<0 \text{ for } 0 < r < \tilde r_2
           \text{ and  } \eta''(r)>0 \text{ for } r > \tilde r_2, \\
\label{da0} &\hat \eta(r)<0 \text{ for }  0 < r < a_0
           \text{ and } \hat\eta(r)>0 \text{ for } r > a_0,\\
\label{d3} &\hat \eta'(r)>0 \text{ for }  0 < r < a_1
           \text{ and } \hat\eta'(r)<0 \text{ for } r > a_1,\\
\label{a4} &0 < a_0 < \tilde r_1 < \tilde r_2 < 2 a_0, \\
\label{a5} &a_0 < a_1 .
\end{align}

We note that it follows from \eqref{da0} that $F_j^{QCF}(\ab)=0$ for $j=-N,\dots,N+1$ and
$\hat{\mathbf\Psi}^F (\ab)=\mathbf 0$  for $\ab=(a_0,\dots,a_0)\in\Real^{2N+1}.$

The local quasicontinuum approximation is simply the case where all of the representative
atoms are continuum,
giving a decoupled system
\begin{equation*}
\hat{\mathbf\Psi}^L(\rb) =\Pb,
\end{equation*}
where $\hat{\mathbf\Psi}^L(\rb) = (\hat{\psi}_{-N}^L(\rb),\dots,\hat{\psi}_{N}^L(\rb))$ is defined by
\begin{equation}\label{L}
\hat{\psi}_i^L(\rb) = \hat\eta(r_i), \qquad i = -N,\dots,N.
\end{equation}
We thus have that
\[
D_j \hat{\psi}^L_i(\rb)=\hat\eta'(r_i)\delta_{ij},
\]
and we see that the local quasicontinuum approximation is unstable when
$r_i>a_1$ for some $i= -N,\dots,N.$

\begin{figure}
\includegraphics[height=2.25in,width=\textwidth]{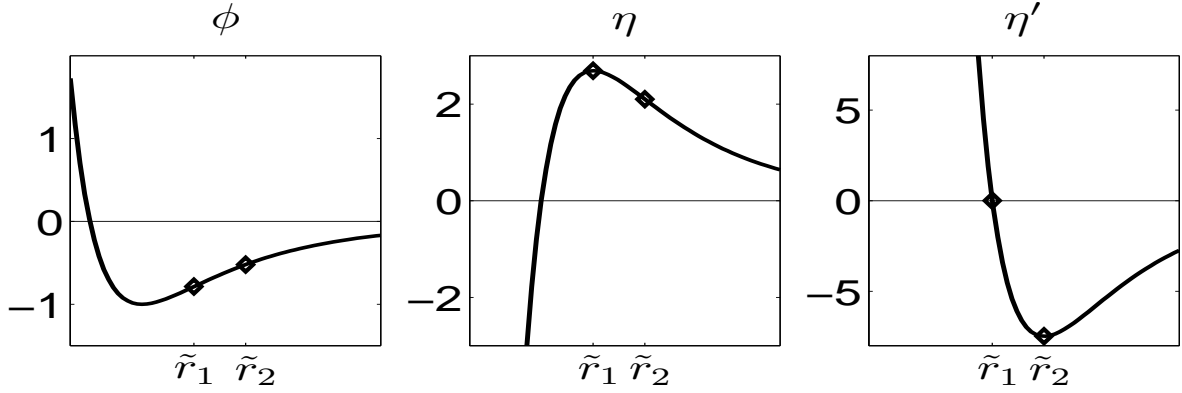}
\caption{\label{phiplot}The Lennard-Jones potential \eqref{lj} demonstrates the
prototypical behavior of $\phi(r)$ and its derivatives, $\eta(r)=\phi'(r)$
and $\eta'(r)=\phi''(r).$}
\end{figure}

\begin{figure}
\includegraphics[height=2.25in,width=\textwidth]{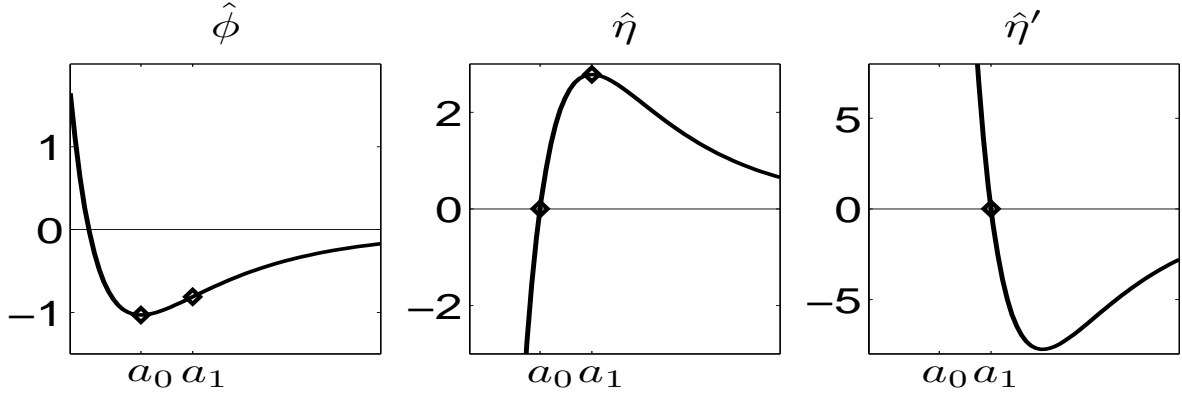}
\caption{\label{hatphiplot}The energy density,
 $\hat \phi(r),$ for the Lennard-Jones potential \eqref{lj} and its derivatives,
 $\hat\eta(r)=\hat\phi'(r)$
and $\hat\eta'(r)=\hat\phi''(r).$}
\end{figure}

\subsection{Existence and Uniqueness by the Inverse Function Theorem}
While the main theorem of this section gives explicit conditions
on $\Pb$ for which $\hat{\mathbf\Psi}^F(\rb)=\Pb$ is solvable,
we begin by showing that the inverse function theorem~\cite{brown} can
be used to show that $\hat{\mathbf\Psi}^F(\rb)$ is bijective in a neighborhood of
$\rb = \ab.$
Therefore, we must analyze
the invertibility of
\begin{equation}
\label{dpsi}
D_j \hat{\psi}^F_i(\rb) = \begin{cases}
    [\eta'(r_i) + 4 \eta'(2 r_i)] \delta_{ij}, &
        -N \leq i \leq -K, \\
   [\eta'(r_i) + \eta'(r_i+r_{i-1}) + \eta'(r_i+r_{i+1})] \delta_{ij}&\\
      \quad +\eta'(r_i + r_{i-1}) \delta_{i-1 j} + \eta'(r_i+r_{i+1})
        \delta_{i+1 j}& \\
    \quad+ [4\eta'(2r_K) - \eta'(r_K+r_{K-1}) - \eta'(r_K+r_{K+1})]
        \delta_{Kj}& \\
    \quad - \eta'(r_K+r_{K-1}) \delta_{K-1j}
        - \eta'(r_K+r_{K+1}) \delta_{K+1j},  & -K+1 \leq i \leq K-1,\\
    [\eta'(r_i) + 4 \eta'(2 r_i)] \delta_{ij}, & K \leq i \leq N.
\end{cases}\end{equation}

\begin{lemma}
\label{local}
If $\eta'(a_0) + 8 \eta'(2 a_0) > 0,$ then $\hat{\mathbf\Psi}^F(\rb)$ is bijective
 in a neighborhood of $\ab.$
\end{lemma}

\begin{proof}
We will show that $D \hat{\mathbf\Psi}^F (\ab)$ (where
$(D \hat{\mathbf\Psi}^F)_{ij} = D_j \hat{\psi}^F_i)$
is nonsingular by demonstrating that it is strictly diagonal
dominant \cite{serre_matrices} (with positive
diagonal), that is,
\begin{equation*}
D_i \hat{\psi}^F_i(\ab) > \sum_{j \neq i} | D_j \hat{\psi}^F_i(\ab)|, \qquad i = -N,\dots,N.
\end{equation*}
Looking at the rows of $D \hat{\mathbf\Psi}^F(\ab),$
all of the nearest-neighbor terms are on the diagonal.  If we collect
all of the next-neighbor terms on the right hand side of the above inequality,
we find that it is sufficient to show that $\eta'(a_0) > 4 |\eta'(2 a_0)|$ for
rows $i=-N,\dots,-K$ and $i=K,\dots,N$ and
$\eta'(a_0) > 8 |\eta'(2 a_0)|$ for rows $i=-K+1,\dots,K-1.$
  Since $2 a_0 > \tilde{r}_1$ by assumption~\eqref{a4},
we have that $\eta'(2a_0)<0$ by~\eqref{d1}.
Thus, the condition
$\eta'(a_0) + 8 \eta'(2 a_0) >0$ implies the strict diagonal dominance
of~\eqref{dpsi}.  Therefore, we have that $D \hat{\mathbf\Psi}^F (\ab)$ is nonsingular,
and so by the inverse function theorem there exists a neighborhood
of $\ab$ in which $\hat{\mathbf\Psi}^F(\rb)$ is bijective~\cite{antman_book}.
\end{proof}

\subsection{Existence and Uniqueness for General External Forces}

In the theorems of this section and the next, we will
use the continuation
method~\cite{antman_book} to find explicit conditions under
which $\hat{\mathbf\Psi}^F(\rb)=\Pb$ has a
unique solution.
The idea of the continuation method is to start with the uncoupled
local quasicontinuum system, construct a homotopy transformation
from the local quasicontinuum system to the
force-based quasicontinuum system, and
show that existence and uniqueness persists through the transformation.
We use the following well-known lemma to show that unique solvability of
the local quasicontinuum system implies existence and uniqueness of a solution to the
force-based quasicontinuum system~\cite{brown}.

\begin{lemma}
\label{homotopy}
Let $\Omega \subset \Real^{2N+1}$ be an open, bounded set.  Suppose that
$\fb, \gb \in C^1(\bar{\Omega};\Real^{2N+1})$ and the homotopy
$\hb(\rb,t) = (1-t) \fb(\rb)+ t \gb(\rb)$ for $t\in[0,1]$ satisfies
\begin{enumerate}
\item $\fb(\rb) = 0$ has
a unique solution in $\Omega,$
\item $\det D_\rb \hb (\rb,t)\neq 0$ in
$\Omega \times [0, 1],$
\item $\hb(\rb,t)\ne 0$
for every $(\rb,t) \in \partial \Omega \times [0,1].$
\end{enumerate}
Then there is a unique solution $\rb \in \Omega$ satisfying $\gb(\rb) = 0.$
\end{lemma}

\begin{proof}
Since  $\det D_\rb \hb(\rb,t) \neq 0$ for all $(\rb,t) \in \Omega\times[0,1]$,
we know by the implicit function theorem~\cite{brown}
 that the solution of $h(\rb,t) = 0$
in a neighborhood of a particular solution $(\rb_*,t_*)$
satisfying $h(\rb_*,t_*) = 0$ can be written as a function $\rb(t)$ in a neighborhood
of $(\rb_*,t_*).$  Thus, by compactness, the curve can be continued until
it leaves the region $\bar \Omega \times [0, 1].$
Since we have assumed that $\hb(\rb,t) \neq 0$ on
$\partial \Omega \times [0,1]$,
there is a one-to-one correspondence between solutions of $\fb(\rb) = 0$
and solutions of $\gb(\rb)=0.$
Therefore, the existence of a
unique solution
of $\fb(\rb) = 0$ for $\rb\in\Omega$ implies
that $\gb(\rb) = 0$ has a unique solution for
$\rb\in\Omega.$
\end{proof}

We define the homotopy
\begin{equation*}
\hb(\rb,t) = (1-t) [\hat{\mathbf\Psi}^L(\rb) - \Pb]
+ t [\hat{\mathbf\Psi}^F(\rb) - \Pb],\qquad t\in[0,1],
\end{equation*}
where $\hat{\mathbf\Psi}^L(\rb)$ is the local quasicontinuum system~\eqref{L}.
The next lemma uses the properties of the local quasicontinuum system
to simplify the hypothesis for the previous lemma.

\begin{lemma}
\label{lem:psi}
Let $r_L$ and $r_U$ satisfy $0<r_L<r_U<a_1,$ and set
$\Omega = (r_L, r_U)^{2N+1}.$
Suppose that for $i=-N,\dots,N$ we have that
\begin{equation}\label{h}
\begin{split}
&h_i(\rb,t) > 0\text{ if }r_i = r_U\text{ and }\rb\in\partial\Omega,\,t\in[0,1],\\
&h_i(\rb,t) < 0\text{ if }r_i = r_L\text{ and }\rb\in\partial\Omega,\,t\in[0,1].
\end{split}
\end{equation}
If $D \hat{\mathbf\Psi}^F(\rb)$ is strictly diagonally dominant
for all $\rb\in\Omega,$ then
there exists a unique solution to $\hat{\mathbf\Psi}^F(\rb) = \Pb.$
\end{lemma}

\begin{proof}
We show that the above conditions are enough to satisfy the hypotheses
of Lemma~\ref{homotopy}.
We first show that $\hat{\mathbf\Psi}^L(\rb) = \Pb$ has a unique solution in $\Omega.$
Since the system is decoupled, we need only demonstrate that the
scalar equations
\begin{equation}\label{de}
\eta(r_j) + 2\eta(2 r_j) = \Phi_j,\qquad j=-N,\dots,N,
\end{equation}
have a unique solution for $r_L<r_j<r_U.$
We have from the hypothesis \eqref{h} on $h_i(\rb,t)$ at $t=0$ that
$\eta(r_L) + 2\eta(2 r_L) - \Phi_j < 0$ and
$\eta(r_U) + 2\eta(2 r_U) - \Phi_j > 0.$  Hence, a solution $r_j$
to~\eqref{de} exists by the intermediate value theorem.
Since
$\hat\eta(r)=\eta(r)+2\eta(2 r)$ is increasing
for $0<r<r_U<a_1$ by~\eqref{d3},
the solution to~\eqref{de}
must be unique.

The hypothesis \eqref{h} implies that
$\hb(\rb,t)\ne 0$
for every $(\rb,t) \in \partial \Omega \times [0,1].$
  Finally, to
show that $\det D_\rb \hb (\rb,t)\neq 0$ for all $(\rb,t)\in
\Omega \times [0, 1]$
it is sufficient to demonstrate
the strict diagonal dominance of $D_\rb \hb(\rb,t)$
for all $(\rb,t)\in
\Omega \times [0, 1].$

Now
\[
D_j h_i(\rb,t) = (1-t) D_j \hat{\psi}^L_i(\rb) + t D_j \hat{\psi}^F_i(\rb),
\]
and
\[
D_j \hat{\psi}^L_i(\rb)=\hat\eta'(r_i)\delta_{ij}.
\]
Since $0<r_i<a_1,$ we have by \eqref{d3} that
$\hat\eta'(r_i) > 0.$  Thus, we have that
$D_j \hat{\psi}^L_i(\rb)$ is strictly diagonal dominant
(with positive diagonal).  We can then conclude that
$D_j h_i(\rb,t)$ is strictly diagonal dominant
(with positive diagonal) if $D_j \hat{\psi}^F_i(\rb)$
is strictly diagonal dominant
(with positive diagonal) since $D_j h_i(\rb,t)$ is then
the sum of two strictly diagonal dominant matrices
(with positive diagonal).
Thus, by
Lemma~\ref{homotopy}, the equation
$\hat{\mathbf\Psi}^F(\rb) - \Pb = 0$ has a unique solution in $\Omega.$
\end{proof}

We next turn to giving results that allow the calculation of
an explicit neighborhood of $\az$ for which $\hat{\mathbf\Psi}^F(\rb)$ is bijective.
  We first give a condition on $0<r_L<r_U$
such that $D \hat{\mathbf\Psi}^F(\rb)$ is strictly diagonally dominant
for all $\rb\in\Omega$ where
$\Omega = (r_L, r_U)^{2N+1}.$

\begin{lemma}\label{reg}
Suppose that $r_L$ and $r_U$ satisfy
\begin{align}
&\frac{\tilde r_2}2<r_L<r_U,\notag\\
&\eta'(r_U)+12\eta'(2r_L) \geq 0.\label{eta}
\end{align}
Then $r_U < a_1,$ and $D \hat{\mathbf\Psi}^F(\rb)$ is strictly
diagonally dominant
for all $\rb\in\Omega$ where
$\Omega = (r_L, r_U)^{2N+1}.$
\end{lemma}
\begin{proof}
First, we note that since $r_L > \frac{\tilde r_2}2,$
\eqref{d2} implies that $12 \eta'(2 r_L) < 12 \eta'(2 r_U).$
Also, \eqref{d2} gives  that the next-neighbor terms are negative,
so that $\eta'(r_U) + 4 \eta'(2 r_U) >
\eta'(r_U) + 12 \eta'(2 r_U) \geq 0.$  Therefore, by \eqref{a4}
we can conclude that $r_U < a_1.$

We can obtain from \eqref{dpsi} by summing all of the
next-nearest-neighbor terms in each row that $D \hat{\mathbf\Psi}^F(\rb)$
is strictly diagonally dominant if
\begin{equation*}
\eta'(r) > 12 |\eta'(2s)| \quad \text{ for all } r,s \in (r_L,r_U).
\end{equation*}
(We note that the hypothesis of Lemma~\ref{local}
required the weaker condition $\eta'(a_0) > 8 |\eta'(2a_0)|$
since we were able to utilize the
cancellation of terms in the expression
\[
D_K \hat{\psi}^F_i(\rb) = 4\eta'(2r_K) - \eta'(r_K+r_{K-1}) - \eta'(r_K+r_{K+1}),
\qquad i = -K+3,\dots,K-3,
\]
when it was evaluated at $\rb=\ab$  to obtain that
 $D_K \hat{\psi}^F_i(\ab)=2\eta'(2a_0)$ for $i = -K+3,\dots,K-3.$)

We have by \eqref{d1} and \eqref{d2} that
\begin{equation}\label{r1}
\eta'(r)>0 \text{ and $\eta'(r)$ is decreasing for }r<\tilde r_1,
\end{equation}
and we have that
\begin{equation}\label{r2}
\eta'(2s)<0 \text{ and $\eta'(2s)$ is increasing for }s>\frac{\tilde r_2}2.
\end{equation}
It follows from \eqref{r2} that to prove
strict diagonal dominance it is sufficient to show that
\begin{equation*}
\eta'(r) + 12 \eta'(2s) > 0 \quad \text{ for all } r,s \in (r_L,r_U).
\end{equation*}
We have from \eqref{r1} and \eqref{r2} that the above
condition follows from the hypothesis \eqref{eta}.
\end{proof}

\begin{theorem}
\label{thm:wp1}
Suppose that $r_L$ and $r_U$ satisfy
\begin{align}
&\frac{\tilde r_2}2<r_L<r_U, \label{what}\\
&\eta'(r_U)+12\eta'(2r_L)\geq0.\label{etaz}
\end{align}
If
\begin{equation*}
\eta(r_L) + 4 \eta(2 r_L) - 2 \eta(2 r_U) < \Phi_j < \eta(r_U)
     + 4 \eta(2 r_U) - 2 \eta(2 r_L), \qquad j=-N,\dots,N,
\end{equation*}
then $\hat{\mathbf\Psi}^F(\rb) = \Pb$
has a unique solution $\rb$ in $\Omega = (r_L, r_U)^{2N+1}.$
\end{theorem}
\begin{proof}
We have by Lemma~\ref{reg} that
$D \hat{\mathbf\Psi}^F(\rb)$ is strictly diagonally dominant
for all $\rb\in\Omega.$
By Lemma~\ref{homotopy}, we need only show
that
\begin{equation*}
\begin{split}
&h_i(\rb,t) > 0\text{ if }r_i = r_U\text{ and }\rb\in\partial\Omega,\,t\in[0,1],\\
&h_i(\rb,t) < 0\text{ if }r_i = r_L\text{ and }\rb\in\partial\Omega,\,t\in[0,1].
\end{split}
\end{equation*}

If we look at the rows of $\hat{\mathbf\Psi}^F(\rb)$, we see that
there is always one nearest-neighbor term and at most four
 positive and two negative next-nearest-neighbor terms~\eqref{stress2}.
 We also note that $2r_L>\tilde r_2>\tilde r_1$ by \eqref{what} and \eqref{a4},
so we have by \eqref{d1} that $\eta(2r)$ is deceasing for $r\ge r_L$ and
\[
2\eta(2 r_L) - 2 \eta(2 r_U)>0.
\]
We can thus estimate $\hb(\rb,t)$ on the boundary to get
\begin{equation}\label{ditz}
\begin{split}
\min_{\rb \in \partial \Omega , r_j = r_U} h_j(\rb,t) &=
    \min_{\rb \in \partial \Omega , r_j = r_U}
     (1-t) [\hat{\psi}_j^L(\rb) - \Phi_j] + t[\hat{\psi}^F_j(\rb) - \Phi_j]
     \\
&\ge (1-t)\left[\eta(r_U)+2\eta(2r_U)-\Phi_j\right]
+t\left[\eta(r_U)+4\eta(2r_U)-2\eta(2r_L)-\Phi_j\right]\\
    & = (1-t) [2 \eta(2 r_L) - 2 \eta(2 r_U)]
    +\left[\eta(r_U)
     + 4 \eta(2 r_U) - 2 \eta(2 r_L)\right]-\Phi_j> 0,\\
\max_{\rb \in \partial \Omega , r_j = r_L} h_j(\rb,t) &=
    \max_{\rb \in \partial \Omega , r_j = r_L}
     (1-t) [\hat{\psi}_j^L(\rb) - \Phi_j] + t[\hat{\psi}^F_j(\rb) - \Phi_j] \\
     &\le (1-t)\left[\eta(r_L)+2\eta(2r_L)-\Phi_j\right]
+t\left[\eta(r_L)+4\eta(2r_L)-2\eta(2r_U)-\Phi_j\right]\\
    & = (1-t) [2 \eta(2 r_U) - 2 \eta(2 r_L)]
    +\left[\eta(r_L) + 4 \eta(2 r_L) - 2 \eta(2 r_U) \right] -\Phi_j< 0.
\end{split}
\end{equation}

Therefore, by Lemma~\ref{lem:psi}, we have a unique solution to
$\hat{\mathbf\Psi}^F(\rb) = \Pb$ in $\Omega.$
\end{proof}

We now consider the case in which
the external
forces, $f_{j},$ are anti-symmetric \eqref{forcesym}, so that
the
external potential $\Phi_j$ is symmetric~\eqref{Phisym}.
We then have that the solution $\rb$ to
$\hat{\mathbf\Psi}^F(\rb) = \Pb$ obtained in the
proof of Theorem~\ref{thm:wp1} is symmetric, and will
also solve the equilibrium equations \eqref{equil}.

\begin{corollary}
\label{cor:wp1}
Suppose that $r_L$ and $r_U$ satisfy
\begin{align*}
&\frac{\tilde r_2}2<r_L<r_U,\\
&\eta'(r_U)+12\eta'(2r_L)\geq0,
\end{align*}
and that
the external
forces, $f_{j},$ are anti-symmetric \eqref{forcesym}.
If
\begin{equation*}
\eta(r_L) + 4 \eta(2 r_L) - 2 \eta(2 r_U) < \Phi_j < \eta(r_U)
     + 4 \eta(2 r_U) - 2 \eta(2 r_L), \qquad j=-N,\dots,N,
\end{equation*}
then the equilibrium equations
\begin{equation*}
F_{j}^{QCF}(\rb)+f_{j}=0,\qquad j=-N,\dots,N+1,
\end{equation*}
have a unique symmetric solution $\rb$ in $\Omega = (r_L, r_U)^{2N+1}.$
\end{corollary}

\begin{proof}
We first note that since the external forces, $f_{j},$ are anti-symmetric,
the external potential $\Pb$ is symmetric.
We consider the solution $\rb(t)$ to
 the homotopy continuation
\begin{equation}\label{d}
\hb\left(\rb(t),t\right) = (1-t) [\hat{\mathbf\Psi}^L\left(\rb(t)\right) - \Pb]
+ t [\hat{\mathbf\Psi}^F\left(\rb(t)\right) - \Pb],\qquad t\in[0,1],
\end{equation}
that we have from Lemma \ref{homotopy}.
The unique solution to
the decoupled local quasicontinuum system
 $\hat{\mathbf\Psi}^L(\rb) =\Pb$ is symmetric, therefore
$\rb(0)$ is symmetric.
We obtain by differentiating \eqref{d} that
\[
 \nabla_{\rb}\hb\left(\rb(t),t\right)\rb_t\left(t\right)
+\hb_t\left(\rb(t),t\right)=0,\qquad t\in[0,1],
\]
where
\[
 \nabla_{\rb}\hb\left(\rb(t),t\right)=(1-t)\nabla_{\rb}\hat{\mathbf\Psi}^L\left(\rb(t)\right)
+t\nabla_{\rb}\hat{\mathbf\Psi}^F\left(\rb(t)\right)
\]
and
\[
 \hb_t\left(\rb(t),t\right)=\hat{\mathbf\Psi}^F\left(\rb(t)\right)-
\hat{\mathbf\Psi}^L\left(\rb(t)\right).
\]
All terms above are symmetric whenever $\rb(t)$ is, and therefore
the symmetry of $\rb(0)$ implies that $\rb_t(t)$ and $\rb(t)$ are
symmetric in the whole interval $[0,1],$ and in particular
$\rb(1)$ is symmetric. The proof that $\rb$ is symmetric
is completed by observing that $\rb(1)$ is the
unique solution to $\hat{\mathbf\Psi}^F(\rb) = \Pb.$

Since the
 unique solution $\rb$ to  $\hat{\mathbf\Psi}^F(\rb) = \Pb$ is symmetric,
 we can conclude that $\rb$ is the unique symmetric solution to
 the force-based equilibrium equations
 \begin{equation*}
F_{j}^{QCF}(\rb)+f_{j}=0,\qquad j=-N,\dots,N+1. \qedhere
\end{equation*}
\end{proof}

We now apply Corollary \ref{cor:wp1} to the Lennard-Jones potential \eqref{lj}.
\begin{corollary}
\label{ljregion}
We assume that
the external
forces, $f_{j},$ are anti-symmetric \eqref{forcesym}.
For any $\frac{\tilde r_2}{2} < r_U < a_1$, let
\begin{equation}\label{ljrl}
r_L = \max \left( \frac{\tilde r_2}{2}, \left( \frac{63}{16 \eta'(r_U)}
        \right)^{\frac{1}{8}} \right).
\end{equation}
If $r_L < r_U$, then
the equilibrium equations
\begin{equation*}
F_{j}^{QCF}(\rb)+f_{j}=0,\qquad j=-N,\dots,N+1,
\end{equation*}
have a unique symmetric solution $\rb$ in $\Omega = (r_L, r_U)^{2N+1}$
whenever
\begin{equation*}
\eta(r_L) + 4 \eta(2 r_L) - 2 \eta(2 r_U) < \Phi_j < \eta(r_U)
     + 4 \eta(2 r_U) - 2 \eta(2 r_L), \qquad j=-N,\dots,N.
\end{equation*}
\end{corollary}

\begin{proof}
To prove this corollary using
 Theorem \ref{thm:wp1}, we need to show that \eqref{etaz} holds.
Evaluating \eqref{etaz} for the Lennard-Jones potential~ \eqref{lj}, we get
\begin{equation*}
\begin{split}
\eta'(r_U) + 12 \eta'(2 r_L) &= \eta'(r_U) + 12 \left[
\frac{156}{2^{14} r_L^{14}} - \frac{84}{2^8 r_L^8} \right] \\
&> \eta'(r_U) - \frac{63}{16 r_L^8}.
\end{split}
\end{equation*}
Therefore, \eqref{etaz} holds if
$r_L \geq \left( \frac{63}{16 \eta'(r_U)}\right)^{1/8}.$
\end{proof}

After solving for $r_L$ in \eqref{ljrl}, we have to check the
additional hypothesis~\eqref{what} that $r_L < r_U$, as it is not true for
every $r_U \in (\hat r_2/2, a_1).$  However, this is not a very
restrictive assumption since it is true for all $r_U < 1.1003$,
whereas $a_1 = 1.1059.$ The lower end for the interval is
$\frac{\tilde r_2}{2} = 0.6085.$

Using Corollary \ref{ljregion}, we now find a symmetric region about
$\Pb = 0.$  Solving numerically, if we take $r_L = 0.9700$ and $r_U = 1.0883$,
we can conclude that for any $\Pb$ satisfying $-2.62 < \Phi < 2.62$, we
can uniquely solve the force-based quasicontinuum equilibrium equations~
\eqref{equil}.  For the Lennard-Jones potential $\hat \eta(a_1) = 2.781$,
so this symmetric region extends quite close to the load limit.

\begin{remark}
The techniques of this section can be applied to the analysis of the
fully atomistic model or the constrained atomistic quasicontinuum
approximation, and the analysis will be simplified as there are no
non-local conjugate forces.  In both cases, the continuation
from the local quasicontinuum approximation can be used.

We also note that the interfacial terms satisfy the convexity
condition
\[
\Su_j(r_{j-1},r_j) = - \half\phi(2r_{j-1}) +  \phi(r_{j-1} + r_j)
    - \half\phi(2r_j)>0,\qquad j=-N+1,\dots,N,
\]
for $\rb$ in the region defined by \eqref{what} in
Theorem~\ref{thm:wp1} and Corollary ~\ref{cor:wp1} since \eqref{r2} holds in this
region.
\end{remark}

\section{Convergence of the Ghost Force Iteration}
\label{iterza}

We now give a similar analysis for the iterative equations
\eqref{iterPsih} which in vector form are
\begin{equation}
\label{iterPsihz}
\hat\pb^E(\rb^{n+1}) + \hat\pb^G(\rb^{n}) = \Pb,
\end{equation}
where $\hat{\mathbf\Psi}^E(\rb)=
\left(\hat\psi^E_{-N}(\rb),\dots,\hat\psi_N^E(\rb)\right)\in\Real^{2N+1}$
has the symmetric form given in \eqref{psiE}, and
$\hat{\mathbf\Psi}^G(\rb)=
\left(\hat\psi^G_{-N}(\rb),\dots,\hat\psi_N^G(\rb)\right)\in\Real^{2N+1}$
has the symmetric form given in \eqref{psiG}.
We also note that we can
ignore $\hat\psi^E_{-N-1},\,\hat\psi^E_{N+1},\,
\hat\psi^G_{-N-1},\,\hat\psi^G_{N+1}$ and
$\Phi_{-N-1},\,\Phi_{N+1}$
in the above formulation since
\[
\hat\psi^E_{-N-1}=\hat\psi^E_{N+1}=\hat\psi^G_{-N-1}=\hat\psi^G_{N+1}=\Phi_{-N-1}=\Phi_{N+1}=0.
\]
We will determine $\Omega \subset \Real^{2N+1}$ and
$D \subset \Real^{2N+1}$ such that there is a unique $\rb^{n+1} \in \Omega$ that
satisfies~\eqref{iterPsihz} whenever
$\rb^n \in \Omega \text{
and } \Pb \in D,$ and we will show that the induced mapping from $\rb^n$ to
$\rb^{n+1}$ is a contraction mapping.

First, we need to compute $D \hat{\mathbf\Psi}^E(\rb)$ and $D \hat{\mathbf\Psi}^G(\rb).$ The
following expressions are rather complex, but all we will need
to recognize from them is that the number of first and second neighbor terms in
each row.  We have
\begin{equation}
\label{eq:Dpsi_E}
D_j \hat\psi^E_i(\rb) =
\begin{cases}
[\eta'(r_i) + 4 \eta'(2r_i)] \delta_{ij}
    ,& K+2 \leq i \leq N  ,\\
[\eta'(r_{i}) + 4\eta(2r_{i}) + \half \eta'(r_i+r_{i-1})]
        \delta_{ij}& \\
    \quad + [\half \eta'(r_i+r_{i-1})] \delta_{i-1j} ,& i = K+1,\\
[\eta'(r_i) + 2 \eta'(2r_i) + \half \eta'(r_i+r_{i-1})
    + \half \eta'(r_i+r_{i+1})] \delta_{ij}& \\
    \quad + [\half \eta'(r_i+r_{i-1})] \delta_{i-1j}
        + [\half \eta'(r_i+r_{i+1})] \delta_{i+1j} ,& i = K ,\\
[\eta'(r_{i}) + \eta'(r_i+r_{i-1}) + \half \eta'(r_i+r_{i+1})]
        \delta_{ij} & \\
    \quad + [\eta'(r_i+r_{i-1})] \delta_{i-1j}
        + [\half \eta'(r_i+r_{i+1})] \delta_{i+1j} ,& i = K-1 ,\\
[\eta'(r_i) + \eta'(r_i+r_{i-1}) + \eta'(r_i+r_{i+1})] \delta_{ij}&\\
    \quad +\eta'(r_i + r_{i-1}) \delta_{i-1 j} + \eta'(r_i+r_{i+1})
        \delta_{i+1 j} ,& -K+2 \leq i \leq K-2,\\
\eta'(r_i)\delta_{ij}+\dots ,& -N\leq i\leq -K+3,
\end{cases}
\end{equation}
and
\begin{equation}
\label{eq:Dpsi_G}
D_j \hat\psi^G_i(\rb) =
\begin{cases}
0 ,& K+2 \leq i \leq N ,\\
- \half \eta'(r_K + r_{K+1}) \delta_{Kj} - \half \eta'(r_K + r_{K+1})
        \delta_{K+1j} ,& i = K+1 ,\\
[2 \eta'(2r_K) - \half \eta'(r_K+r_{K-1}) - \half \eta'(r_K+r_{K+1})]
        \delta_{Kj}&\\
    \quad - \half \eta'(r_K+r_{K-1}) \delta_{K-1j} - \half \eta'(r_K+r_{K+1})
        \delta_{K+1j} ,& i = K ,\\
[4 \eta'(2r_K) - \half \eta'(r_K+r_{K-1}) - \eta'(r_K+r_{K+1})] \delta_{Kj}&\\
    \quad - \half \eta'(r_K+r_{K-1}) \delta_{K-1j} - \eta'(r_K+r_{K+1})
        \delta_{K+1j} ,& i = K-1 ,\\
[4 \eta'(2r_K) - \eta'(r_K+r_{K-1}) - \eta'(r_K+r_{K+1})] \delta_{Kj}&\\
    \quad - \eta'(r_K+r_{K-1}) \delta_{K-1j} - \eta'(r_K+r_{K+1})
        \delta_{K+1j} ,& -K+2 \leq i \leq K-2,\\
\dots ,& -N\leq i\leq -K+3.
\end{cases}
\end{equation}

We recall that the maximum norm for $\rb\in\Real^{2N+1}$ is given by
\[
\|\rb\|_\infty =
\max_{i=-N,\dots,N} |r_i|.
\]

\begin{theorem}
\label{thm:wp2}
Suppose that $r_L$ and $r_U$ satisfy
\begin{align}
&\frac{\hat r_2}{2} <  r_L < r_U,\notag \\
&\label{thirteen} \eta'(r_U) + 13 \eta'(2 r_L) > 0,
\end{align}
and that $\Pb$ satisfies
\begin{equation} \label{eq:phirange2}
\eta(r_L) + 4 \eta(2 r_L) - 2 \eta(2 r_U) < \Phi_j < \eta(r_U)
     + 4 \eta(2 r_U) - 2 \eta(2 r_L), \qquad j = -N,\dots,N.
\end{equation}
Then for every $\rb^{n} \in \Omega = (r_L,r_U)^{2N+1}$
there is a unique $\rb^{n+1} \in \Omega$
such that
\begin{equation}\label{it}
\hat{\mathbf\Psi}^E(\rb^{n+1}) + \hat{\mathbf\Psi}^G(\rb^{n}) = \Pb.
\end{equation}
We also have that the
induced mapping $\rb^n \rightarrow \rb^{n+1}$ is a contraction and
satisfies the inequality
\begin{equation} \label{contract}
\maxnorm{\rb^{n+1} - \sbf^{n+1}}\le \frac{  8 |\eta'(2 r_L)|}{\eta'(r_U) -  5 |\eta'(2 r_L)|}
\maxnorm{\rb^n - \sbf^n},
\end{equation}
where we have from \eqref{thirteen} that
\begin{equation*}
\frac{  8 |\eta'(2 r_L)|}{\eta'(r_U) -  5 |\eta'(2 r_L)|} < 1.
\end{equation*}
\end{theorem}

\begin{proof}
The first part of the proof will be very similar to the proofs of
Lemma~\ref{reg} and Theorem~\ref{thm:wp1}, as we will again satisfy
the hypotheses of
Lemma~\ref{lem:psi}.  First, we note that \eqref{thirteen} implies
that $r_U < a_1$ by Lemma \ref{reg}.
We prove strict diagonal dominance for $D_j \hat\psi^E_i(\rb)$
given by \eqref{eq:Dpsi_E}.  For this argument,
we only need to show that $\eta'(r) > 5 |\eta'(2s)|$
whenever $r,s \in (r_L,r_U),$ or since $\eta'(2s) < 0$,
we need only show
\[
\eta'(r) + 5 \eta'(2s) > 0 \quad \text{for } r,s \in (r_L,r_U).
\]
We will need the factor $13$ in the hypothesis
\eqref{thirteen} to prove that the mapping is a contraction, which
is why the hypothesis is as strong as it is.
We further note that from \eqref{r1},
\eqref{r2}, and the hypothesis \eqref{eq:phirange2} that we have
\begin{equation*}
\begin{split}
\eta'(r) + 5 \eta'(2s) &> \eta'(r) + 13 \eta'(2s)\\
                       &> \eta'(r_U) + 13 \eta'(2 r_L)>0,\qquad
r,s \in (r_L,r_U).
\end{split}
\end{equation*}
Thus, we have established that $\det D \hat{\mathbf\Psi}^E(\rb)> 0$ in
$\Omega.$

We now verify that $\hb(\rb,t)$ satisfies the condition ~\eqref{h},
that is, it does not vanish on
$\partial\Omega\times [0,1].$
In the proof of Theorem ~\ref{thm:wp1}, we analyzed $\hb(\rb,t)$
on $\partial\Omega\times [0,1]$ by grouping
the second-neighbor terms.  From this point of view~\eqref{ditz}, the iteration problem
$\hat{\mathbf\Psi}^E(\rb^{n+1}) + \hat{\mathbf\Psi}^G(\rb^{n})=\Pb,$
is identical to the problem, $\hat{\mathbf\Psi}(\rb) =\Pb.$
Thus, we have that
\begin{equation*}
\begin{split}
\min_{\substack{\rb^{n+1} \in \partial \Omega ,\\ r_j = r_U}} h_j(\rb^{n+1})
&\geq \min_{\substack{\rb^{n+1},\rb^{n} \in \partial \Omega ,\\ r_j = r_U}}
     (1-t) [\hat\psi_j^L(\rb^{n+1}) - \Phi_j] + t[\hat\psi^E_j(\rb^{n+1})
     + \hat\psi^G(\rb^n) - \Phi_j] \\
    & \ge (1-t) [2 \eta(2 r_L) - 2 \eta(2 r_U)]
    +t \left[\eta(r_U)
     + 4 \eta(2 r_U) - 2 \eta(2 r_L)-\Phi_j\right]> 0,\\
\max_{\substack{\rb^{n+1} \in \partial \Omega ,\\ r_j = r_L}} h_j(\rb^{n+1})
&\leq \max_{\substack{\rb^{n+1},\rb^{n} \in \partial \Omega ,\\ r_j = r_L}}
     (1-t) [\hat\psi_j^L(\rb^{n+1}) - \Phi_j] + t[\hat\psi^E_j(\rb^{n+1})
     + \hat\psi^G(\rb^n) - \Phi_j] \\
    & \le (1-t) [2 \eta(2 r_U) - 2 \eta(2 r_L)]
    +t \left[\eta(r_L) + 4 \eta(2 r_L) - 2 \eta(2 r_U)
    -\Phi_j \right]< 0.
\end{split}
\end{equation*}
Therefore, we have proven that $\hb(\rb,t)$
satisfies the condition ~\eqref{h},
that is, it does not vanish on
$\partial\Omega\times [0,1].$
We have now verified the hypotheses of Lemma~\ref{lem:psi}
to conclude the existence of a unique solution $\rb^{n+1}\in\Omega$
to the iteration equation ~\eqref{it}.

To prove that the mapping~\eqref{it} is a contraction,
we suppose
$\rb^n, \rb^{n+1}, \sbf^n, \sbf^{n+1} \in \Omega$ satisfy
\begin{equation*}
\begin{split}
\hat{\mathbf\Psi}^E(\rb^{n+1}) + \hat{\mathbf\Psi}^G(\rb^n) = \Pb, \\
\hat{\mathbf\Psi}^E(\sbf^{n+1}) + \hat{\mathbf\Psi}^G(\sbf^n) = \Pb.
\end{split}
\end{equation*}
We then have that
\begin{equation}\label{contract2}
\hat{\mathbf\Psi}^E(\rb^{n+1}) - \hat{\mathbf\Psi}^E(\sbf^{n+1}) = \hat{\mathbf\Psi}^G(\sbf^n) - \hat{\mathbf\Psi}^G(\rb^n).
\end{equation}

By the fundamental theorem of calculus, we have that
\begin{equation}\label{har}
\begin{split}
\hat{\mathbf\Psi}^E(\rb^{n+1}) - \hat{\mathbf\Psi}^E(\sbf^{n+1})&=L^E(\rb^{n+1},\sbf^{n+1})(\rb^{n+1}-\sbf^{n+1}),\\
\hat{\mathbf\Psi}^G(\rb^{n}) - \hat{\mathbf\Psi}^G(\sbf^{n})&=L^G(\rb^{n},\sbf^{n})(\rb^{n}-\sbf^{n}),
\end{split}
\end{equation}
where
\begin{equation*}
\begin{split}
 L^E(\rb^{n+1},\sbf^{n+1})&=\int_0^1 D\hat{\mathbf\Psi}^E\left(\rb^{n+1} + \theta
 (\sbf^{n+1}-\rb^{n+1})\right)\,d\theta,\\
L^G(\rb^{n},\sbf^{n})&=\int_0^1 D\hat{\mathbf\Psi}^G\left(\rb^{n} + \theta (\sbf^{n}-\rb^{n})\right)\,d\theta.
\end{split}
\end{equation*}
We then have by \eqref{contract2} and \eqref{har} that
\begin{equation*}
L^E(\rb^{n+1},\sbf^{n+1})(\rb^{n+1}-\sbf^{n+1})=-L^G(\rb^{n},\sbf^{n})
(\rb^{n}-\sbf^{n}).
\end{equation*}

To show that $L^E(\rb^{n+1},\sbf^{n+1})$ is nonsingular
and estimate its inverse, we define
the diagonal matrix, $L^D(\rb^{n+1},\sbf^{n+1}),$
that contains the dominant nearest-neighbor terms
of $L^E(\rb^{n+1},\sbf^{n+1})$ by $\left(L^D(\rb^{n+1},\sbf^{n+1})\right)_{ij}
=L^D_{ij}(\rb^{n+1},\sbf^{n+1})$ where
\begin{equation*}
 L^D_{ij}(\rb^{n+1},\sbf^{n+1})=\left[\int_0^1 \eta'\left(r_i^{n+1} + \theta
 (s_i^{n+1}-r_i^{n+1})\right)\,d\theta\right]
\delta_{ij}.
\end{equation*}

Since $\rb^{n} + \theta
 (\sbf^{n}-\rb^{n})\in\Omega=(r_L, r_U)^{2N+1}$ if $\rb^{n},\,\sbf^{n}\in\Omega$
and $\theta\in(0,1),$ we can conclude from \eqref{r1} and \eqref{r2} that
\begin{equation}\label{br}
\begin{split}
\maxnorm{L^D(\rb^{n+1},\sbf^{n+1})^{-1}} &\leq \frac{1}{\eta'(r_U)}, \\
\maxnorm{ L^E(\rb^{n+1},\sbf^{n+1})- L^D(\rb^{n+1},\sbf^{n+1}) } &\leq 5 |\eta'(2 r_L)|,
\end{split}
\end{equation}
where the matrix norm that is induced by the $\maxnorm{\rb}$ vector norm is
\[
\|L\|_\infty = \max_{i=-N,\dots,N} \sum_{j=-N}^N |L_{ij}|.
\]
We have that
\[
 L^E(\rb^{n+1},\sbf^{n+1})=L^D(\rb^{n+1},\sbf^{n+1})
\left[I+L^D(\rb^{n+1},\sbf^{n+1})^{-1}\left(L^E(\rb^{n+1},\sbf^{n+1})
-L^D(\rb^{n+1},\sbf^{n+1})\right)\right],
\]
so it follows from \eqref{br} that $L^E(\rb^{n+1},\sbf^{n+1})$ is
nonsingular and we have the estimate
\begin{equation}\label{dar}
 \maxnorm{L^E(\rb^{n+1},\sbf^{n+1})^{-1}}\le
\frac{\maxnorm{ L^D(\rb^{n+1},\sbf^{n+1})^{-1}}}{1-\maxnorm{ L^D(\rb^{n+1},\sbf^{n+1})^{-1}}
\maxnorm{ L^E(\rb^{n+1},\sbf^{n+1})- L^D(\rb^{n+1},\sbf^{n+1}) }}.
\end{equation}
Hence, can state that
\begin{equation}\label{ya}
 \rb^{n+1}-\sbf^{n+1}=-\left[L^E(\rb^{n+1},\sbf^{n+1})\right]^{-1}
L^G(\rb^{n},\sbf^{n})
(\rb^{n}-\sbf^{n}).
\end{equation}

From \eqref{eq:Dpsi_G}, we can obtain from \eqref{r2} the estimate
\begin{equation*}
 \maxnorm{L^G(\rb^{n},\sbf^{n})} \leq 8 |\eta'(2 r_L)|,
\end{equation*}
so we have from \eqref{ya} and \eqref{dar} that
\begin{equation*}
\begin{split}
\maxnorm{\rb^{n+1} - \sbf^{n+1}} &\leq \maxnorm{L^E(\rb^{n+1},\sbf^{n+1})^{-1}L^G(\rb^{n},\sbf^{n})}
    \maxnorm{\rb^n - \sbf^n}\\
&\leq \maxnorm{L^E(\rb^{n+1},\sbf^{n+1})^{-1}}
    \maxnorm{L^G(\rb^{n},\sbf^{n})} \maxnorm{\rb^n - \sbf^n}\\
&\leq \frac{  8 |\eta'(2 r_L)|}{\eta'(r_U) -  5 |\eta'(2 r_L)|}\maxnorm{\rb^n - \sbf^n}.
\end{split}
\end{equation*}
We have from \eqref{thirteen} that
\begin{equation*}
\frac{  8 |\eta'(2 r_L)|}{\eta'(r_U) -  5 |\eta'(2 r_L)|} < 1,
\end{equation*}
so the mapping $\rb^n \rightarrow \rb^{n+1}$ is a contraction and
\[
\maxnorm{\rb^{n+1} - \sbf^{n+1}}\le \frac{  8 |\eta'(2 r_L)|}{\eta'(r_U) -  5 |\eta'(2 r_L)|}
\maxnorm{\rb^n - \sbf^n}. \qedhere
\]
\end{proof}

We can prove the following corollary of Theorem~\ref{thm:wp2} for the ghost force iteration
by an argument similar to that used to derive Corollary~\ref{cor:wp1} from Theorem~\ref{thm:wp1}.

\begin{corollary}
\label{cor:wp2}
Suppose that
the external
forces, $f_{j},$ are anti-symmetric \eqref{forcesym},
that $r_L$ and $r_U$ satisfy
\begin{align}
&\frac{\hat r_2}{2} <  r_L < r_U , \notag \\
&\label{thirteenzz} \eta'(r_U) + 13 \eta'(2 r_L) > 0,
\end{align}
and that $\Pb$ satisfies
\begin{equation*}
\eta(r_L) + 4 \eta(2 r_L) - 2 \eta(2 r_U) < \Phi_j < \eta(r_U)
     + 4 \eta(2 r_U) - 2 \eta(2 r_L), \qquad j = -N,\dots,N.
\end{equation*}

Then for every symmetric $\rb^{n} \in \Omega = (r_L,r_U)^{2N+1}$
there is a unique symmetric $\rb^{n+1} \in \Omega$
such that
\begin{equation*}
F_j^{QCE}(\rb^{n+1}) + F_j^G(\rb^{n}) + f_j = 0, \qquad j = -N,\dots,N+1.
\end{equation*}
We also have that the
induced mapping $\rb^n \rightarrow \rb^{n+1}$ is a contraction and
satisfies the inequality
\begin{equation}\label{contract3}
\maxnorm{\rb^{n+1} - \sbf^{n+1}}\le \frac{  8 |\eta'(2 r_L)|}{\eta'(r_U) -  5 |\eta'(2 r_L)|}
\maxnorm{\rb^n - \sbf^n},
\end{equation}
where we have from \eqref{thirteenzz} that
\begin{equation*}
\frac{  8 |\eta'(2 r_L)|}{\eta'(r_U) -  5 |\eta'(2 r_L)|} < 1.
\end{equation*}
The mapping $\rb^n \rightarrow \rb^{n+1}$ converges to the unique symmetric $\rb$
in $\Omega$ that satisfies the force-based quasicontinuum equations
\begin{equation*}\label{f}
F_{j}^{QCF}(\rb)+f_{j}=0,\qquad j=-N,\dots,N+1.
\end{equation*}
\end{corollary}

We now apply Corollary \ref{cor:wp2} to the Lennard-Jones potential
as in the previous section.
This time, we not only need to satisfy the basic
inequality~\eqref{thirteen} on $r_L$ and $r_U,$
but we also need to verify that the contraction constant
in \eqref{contract3} is less than 1.  So, we pick $\gamma \in (0,1)$ and
solve for
\begin{equation*}
\frac{  8 |\eta'(2 r_L)|}{\eta'(r_U) -  5 |\eta'(2 r_L)|} < \gamma.
\end{equation*}
Using the same argument as in Corollary \ref{ljregion}, for any
$r_U \in(\hat r_2/2, a_1),$ we choose
\[
 r_L = \max \left( \frac{\tilde r_2}{2}, \left(
 \frac{156 (5 + 8/\gamma)}{256 \eta'(r_U)}
        \right)^{\frac{1}{8}} \right).
\]
If the resulting $r_L$ is less than $r_U$, we then have a region for
symmetric $\Pb$ for which the iteration is well-defined and a contraction.

Using the above with contraction constant $\gamma = \frac{1}{2},$ we find that for any symmetric
$\Pb \in (-2.56, 2.56)^{2N+1}$ and
symmetric $\rb^n \in (.9706,1.0771)^{2N+1}$ there is a unique
symmetric $\rb^{n+1}\in (.9706,1.0771)^{2N+1}$ that satisfies the ghost force iteration equations
\begin{equation*}
F_j^{QCE}(\rb^{n+1}) + F_j^G(\rb^{n}) + f_j = 0, \qquad j = -N,\dots,N+1.
\end{equation*}
We can finally obtain by taking $\sbf^n=\rb,$
where $\rb$ is the unique symmetric solution to the force-based quasicontinuum equations,
\[
F_{j}^{QCF}(\rb)+f_{j}=0,\qquad j=-N,\dots,N+1,
\]
that
\begin{equation*}
\maxnorm{\rb^{n+1} - \rb}\le \frac12
\maxnorm{\rb^n - \rb}.
\end{equation*}

\bibliography{qcf}
\bibliographystyle{abbrv}
}

\end{document}